\documentclass[10pt]{article}
\usepackage{amssymb,amsmath,amsfonts,bbm,pifont,upgreek,bbold,accents}  
\usepackage[colorlinks=true]{hyperref}
\hypersetup{urlcolor=blue, citecolor=red}

\usepackage[dvips]{graphicx}

\ifx\fiverm\undefined

	\newfont\fiverm{cmr5}

\fi

\catcode`@=11 \catcode`!=11

\expandafter\ifx\csname fiverm\endcsname\relax
  \let\fiverm\fivrm
\fi
  
\let\!latexendpicture=\endpicture 
\let\!latexframe=\frame
\let\!latexlinethickness=\linethickness
\let\!latexmultiput=\multiput
\let\!latexput=\put
 
\def\@picture(#1,#2)(#3,#4){%
  \@picht #2\unitlength
  \setbox\@picbox\hbox to #1\unitlength\bgroup 
  \let\endpicture=\!latexendpicture
  \let\frame=\!latexframe
  \let\linethickness=\!latexlinethickness
  \let\multiput=\!latexmultiput
  \let\put=\!latexput
  \hskip -#3\unitlength \lower #4\unitlength \hbox\bgroup}

\catcode`@=12 \catcode`!=12

\catcode`!=11 
 
\def\PiC{P\kern-.12em\lower.5ex\hbox{I}\kern-.075emC}
\def\PiCTeX{\PiC\kern-.11em\TeX}

\def\!ifnextchar#1#2#3{%
  \let\!testchar=#1%
  \def\!first{#2}%
  \def\!second{#3}%
  \futurelet\!nextchar\!testnext}
\def\!testnext{%
  \ifx \!nextchar \!spacetoken 
    \let\!next=\!skipspacetestagain
  \else
    \ifx \!nextchar \!testchar
      \let\!next=\!first
    \else 
      \let\!next=\!second 
    \fi 
  \fi
  \!next}
\def\\{\!skipspacetestagain} 
  \expandafter\def\\ {\futurelet\!nextchar\!testnext} 
\def\\{\let\!spacetoken= } \\  

\def\!tfor#1:=#2\do#3{%
  \edef\!fortemp{#2}%
  \ifx\!fortemp\!empty 
    \else
    \!tforloop#2\!nil\!nil\!!#1{#3}%
  \fi}
\def\!tforloop#1#2\!!#3#4{%
  \def#3{#1}%
  \ifx #3\!nnil
    \let\!nextwhile=\!fornoop
  \else
    #4\relax
    \let\!nextwhile=\!tforloop
  \fi 
  \!nextwhile#2\!!#3{#4}}

\def\!etfor#1:=#2\do#3{%
  \def\!!tfor{\!tfor#1:=}%
  \edef\!!!tfor{#2}%
  \expandafter\!!tfor\!!!tfor\do{#3}}

\def\!cfor#1:=#2\do#3{%
  \edef\!fortemp{#2}%
  \ifx\!fortemp\!empty 
  \else
    \!cforloop#2,\!nil,\!nil\!!#1{#3}%
  \fi}
\def\!cforloop#1,#2\!!#3#4{%
  \def#3{#1}%
  \ifx #3\!nnil
    \let\!nextwhile=\!fornoop 
  \else
    #4\relax
    \let\!nextwhile=\!cforloop
  \fi
  \!nextwhile#2\!!#3{#4}}

\def\!ecfor#1:=#2\do#3{%
  \def\!!cfor{\!cfor#1:=}%
  \edef\!!!cfor{#2}%
  \expandafter\!!cfor\!!!cfor\do{#3}}

\def\!empty{}
\def\!nnil{\!nil}
\def\!fornoop#1\!!#2#3{}

\def\!ifempty#1#2#3{%
  \edef\!emptyarg{#1}%
  \ifx\!emptyarg\!empty
    #2%
  \else
    #3%
  \fi}
 
\def\!getnext#1\from#2{%
  \expandafter\!gnext#2\!#1#2}%
\def\!gnext\\#1#2\!#3#4{%
  \def#3{#1}%
  \def#4{#2\\{#1}}%
  \ignorespaces}

\def\!getnextvalueof#1\from#2{%
  \expandafter\!gnextv#2\!#1#2}%
\def\!gnextv\\#1#2\!#3#4{%
  #3=#1%
  \def#4{#2\\{#1}}%
  \ignorespaces}

\def\!copylist#1\to#2{%
  \expandafter\!!copylist#1\!#2}
\def\!!copylist#1\!#2{%
  \def#2{#1}\ignorespaces}

\def\!wlet#1=#2{%
  \let#1=#2 
  \wlog{\string#1=\string#2}}
 
\def\!listaddon#1#2{%
  \expandafter\!!listaddon#2\!{#1}#2}
\def\!!listaddon#1\!#2#3{%
  \def#3{#1\\#2}}
 
\def\!rightappend#1\withCS#2\to#3{\expandafter\!!rightappend#3\!#2{#1}#3}
\def\!!rightappend#1\!#2#3#4{\def#4{#1#2{#3}}}

\def\!leftappend#1\withCS#2\to#3{\expandafter\!!leftappend#3\!#2{#1}#3}
\def\!!leftappend#1\!#2#3#4{\def#4{#2{#3}#1}}

\def\!lop#1\to#2{\expandafter\!!lop#1\!#1#2}
\def\!!lop\\#1#2\!#3#4{\def#4{#1}\def#3{#2}}

\def\!loop#1\repeat{\def\!body{#1}\!iterate}
\def\!iterate{\!body\let\!next=\!iterate\else\let\!next=\relax\fi\!next}
 
\def\!!loop#1\repeat{\def\!!body{#1}\!!iterate}
\def\!!iterate{\!!body\let\!!next=\!!iterate\else\let\!!next=\relax\fi\!!next}
 
\def\!removept#1#2{\edef#2{\expandafter\!!removePT\the#1}}
{\catcode`p=12 \catcode`t=12 \gdef\!!removePT#1pt{#1}}

\def\placevalueinpts of <#1> in #2 {%
  \!removept{#1}{#2}}
 
\def\!mlap#1{\hbox to 0pt{\hss#1\hss}}
\def\!vmlap#1{\vbox to 0pt{\vss#1\vss}}
 
\def\!not#1{%
  #1\relax
    \!switchfalse
  \else
    \!switchtrue
  \fi
  \if!switch
  \ignorespaces}

\let\!!!wlog=\wlog              
\def\wlog#1{}    

\newdimen\headingtoplotskip     
\newdimen\linethickness         
\newdimen\longticklength        
\newdimen\plotsymbolspacing     
\newdimen\shortticklength       
\newdimen\stackleading          
\newdimen\tickstovaluesleading  
\newdimen\totalarclength        
\newdimen\valuestolabelleading  

\newbox\!boxA                   
\newbox\!boxB                   
\newbox\!picbox                 
\newbox\!plotsymbol             
\newbox\!putobject              
\newbox\!shadesymbol            

\newcount\!countA               
\newcount\!countB               
\newcount\!countC               
\newcount\!countD               
\newcount\!countE               
\newcount\!countF               
\newcount\!countG               
\newcount\!fiftypt              
\newcount\!intervalno           
\newcount\!npoints              
\newcount\!nsegments            
\newcount\!ntemp                
\newcount\!parity               
\newcount\!scalefactor          
\newcount\!tfs                  
\newcount\!tickcase             

\newdimen\!Xleft                
\newdimen\!Xright               
\newdimen\!Xsave                
\newdimen\!Ybot                 
\newdimen\!Ysave                
\newdimen\!Ytop                 
\newdimen\!angle                
\newdimen\!arclength            
\newdimen\!areabloc             
\newdimen\!arealloc             
\newdimen\!arearloc             
\newdimen\!areatloc             
\newdimen\!bshrinkage           
\newdimen\!checkbot             
\newdimen\!checkleft            
\newdimen\!checkright           
\newdimen\!checktop             
\newdimen\!dimenA               
\newdimen\!dimenB               
\newdimen\!dimenC               
\newdimen\!dimenD               
\newdimen\!dimenE               
\newdimen\!dimenF               
\newdimen\!dimenG               
\newdimen\!dimenH               
\newdimen\!dimenI               
\newdimen\!distacross           
\newdimen\!downlength           
\newdimen\!dp                   
\newdimen\!dshade               
\newdimen\!dxpos                
\newdimen\!dxprime              
\newdimen\!dypos                
\newdimen\!dyprime              
\newdimen\!ht                   
\newdimen\!leaderlength         
\newdimen\!lshrinkage           
\newdimen\!midarclength         
\newdimen\!offset               
\newdimen\!plotheadingoffset    
\newdimen\!plotsymbolxshift     
\newdimen\!plotsymbolyshift     
\newdimen\!plotxorigin          
\newdimen\!plotyorigin          
\newdimen\!rootten              
\newdimen\!rshrinkage           
\newdimen\!shadesymbolxshift    
\newdimen\!shadesymbolyshift    
\newdimen\!tenAa                
\newdimen\!tenAc                
\newdimen\!tenAe                
\newdimen\!tshrinkage           
\newdimen\!uplength             
\newdimen\!wd                   
\newdimen\!wmax                 
\newdimen\!wmin                 
\newdimen\!xB                   
\newdimen\!xC                   
\newdimen\!xE                   
\newdimen\!xM                   
\newdimen\!xS                   
\newdimen\!xaxislength          
\newdimen\!xdiff                
\newdimen\!xleft                
\newdimen\!xloc                 
\newdimen\!xorigin              
\newdimen\!xpivot               
\newdimen\!xpos                 
\newdimen\!xprime               
\newdimen\!xright               
\newdimen\!xshade               
\newdimen\!xshift               
\newdimen\!xtemp                
\newdimen\!xunit                
\newdimen\!xxE                  
\newdimen\!xxM                  
\newdimen\!xxS                  
\newdimen\!xxloc                
\newdimen\!yB                   
\newdimen\!yC                   
\newdimen\!yE                   
\newdimen\!yM                   
\newdimen\!yS                   
\newdimen\!yaxislength          
\newdimen\!ybot                 
\newdimen\!ydiff                
\newdimen\!yloc                 
\newdimen\!yorigin              
\newdimen\!ypivot               
\newdimen\!ypos                 
\newdimen\!yprime               
\newdimen\!yshade               
\newdimen\!yshift               
\newdimen\!ytemp                
\newdimen\!ytop                 
\newdimen\!yunit                
\newdimen\!yyE                  
\newdimen\!yyM                  
\newdimen\!yyS                  
\newdimen\!yyloc                
\newdimen\!zpt                  

\newif\if!axisvisible           
\newif\if!gridlinestoo          
\newif\if!keepPO                
\newif\if!placeaxislabel        
\newif\if!switch                
\newif\if!xswitch               

\newtoks\!axisLaBeL             
\newtoks\!keywordtoks           

\newwrite\!replotfile           

\newhelp\!keywordhelp{The keyword mentioned in the error message in unknown. 
Replace NEW KEYWORD in the indicated response by the keyword that 
should have been specified.}    

\!wlet\!!origin=\!xM                   
\!wlet\!!unit=\!uplength               
\!wlet\!Lresiduallength=\!dimenG       
\!wlet\!Rresiduallength=\!dimenF       
\!wlet\!axisLength=\!distacross        
\!wlet\!axisend=\!ydiff                
\!wlet\!axisstart=\!xdiff              
\!wlet\!axisxlevel=\!arclength         
\!wlet\!axisylevel=\!downlength        
\!wlet\!beta=\!dimenE                  
\!wlet\!gamma=\!dimenF                 
\!wlet\!shadexorigin=\!plotxorigin     
\!wlet\!shadeyorigin=\!plotyorigin     
\!wlet\!ticklength=\!xS                
\!wlet\!ticklocation=\!xE              
\!wlet\!ticklocationincr=\!yE          
\!wlet\!tickwidth=\!yS                 
\!wlet\!totalleaderlength=\!dimenE     
\!wlet\!xone=\!xprime                  
\!wlet\!xtwo=\!dxprime                 
\!wlet\!ySsave=\!yM                    
\!wlet\!ybB=\!yB                       
\!wlet\!ybC=\!yC                       
\!wlet\!ybE=\!yE                       
\!wlet\!ybM=\!yM                       
\!wlet\!ybS=\!yS                       
\!wlet\!ybpos=\!yyloc                  
\!wlet\!yone=\!yprime                  
\!wlet\!ytB=\!xB                       
\!wlet\!ytC=\!xC                       
\!wlet\!ytE=\!downlength               
\!wlet\!ytM=\!arclength                
\!wlet\!ytS=\!distacross               
\!wlet\!ytpos=\!xxloc                  
\!wlet\!ytwo=\!dyprime                 

\!zpt=0pt                              
\!xunit=1pt
\!yunit=1pt
\!arearloc=\!xunit
\!areatloc=\!yunit
\!dshade=5pt
\!leaderlength=24in
\!tfs=256                              
\!wmax=5.3pt                           
\!wmin=2.7pt                           
\!xaxislength=\!xunit
\!xpivot=\!zpt
\!yaxislength=\!yunit 
\!ypivot=\!zpt
  \!dimenA=50pt \!fiftypt=\!dimenA     

\!rootten=3.162278pt                   
\!tenAa=8.690286pt                     
\!tenAc=2.773839pt                     
\!tenAe=2.543275pt                     

\def\!cosrotationangle{1}      
\def\!sinrotationangle{0}      
\def\!xpivotcoord{0}           
\def\!xref{0}                  
\def\!xshadesave{0}            
\def\!ypivotcoord{0}           
\def\!yref{0}                  
\def\!yshadesave{0}            
\def\!zero{0}                  

\let\wlog=\!!!wlog
\def\normalgraphs{%
  \longticklength=.4\baselineskip
  \shortticklength=.25\baselineskip
  \tickstovaluesleading=.25\baselineskip
  \valuestolabelleading=.8\baselineskip
  \linethickness=.4pt
  \stackleading=.17\baselineskip
  \headingtoplotskip=1.5\baselineskip
  \visibleaxes
  \ticksout
  \nogridlines
  \unloggedticks}
\def\setplotarea x from #1 to #2, y from #3 to #4 {%
  \!arealloc=\!M{#1}\!xunit \advance \!arealloc -\!xorigin
  \!areabloc=\!M{#3}\!yunit \advance \!areabloc -\!yorigin
  \!arearloc=\!M{#2}\!xunit \advance \!arearloc -\!xorigin
  \!areatloc=\!M{#4}\!yunit \advance \!areatloc -\!yorigin
  \!initinboundscheck
  \!xaxislength=\!arearloc  \advance\!xaxislength -\!arealloc
  \!yaxislength=\!areatloc  \advance\!yaxislength -\!areabloc
  \!plotheadingoffset=\!zpt
  \!dimenput {{\setbox0=\hbox{}\wd0=\!xaxislength\ht0=\!yaxislength\box0}}
     [bl] (\!arealloc,\!areabloc)}
\def\visibleaxes{%
  \def\!axisvisibility{\!axisvisibletrue}}

\def\!fixkeyword#1{%
  \errhelp=\!keywordhelp
  \errmessage{Unrecognized keyword `#1': \the\!keywordtoks{NEW KEYWORD}'}}

\!keywordtoks={enter `i\fixkeyword}

\def\fixkeyword#1{%
  \!nextkeyword#1 }

\def\axis {%
  \def\!nextkeyword##1 {%
    \expandafter\ifx\csname !axis##1\endcsname \relax
      \def\!next{\!fixkeyword{##1}}%
    \else
      \def\!next{\csname !axis##1\endcsname}%
    \fi
    \!next}%
  \!offset=\!zpt
  \!axisvisibility
  \!placeaxislabelfalse
  \!nextkeyword}

\def\!axisbottom{%
  \!axisylevel=\!areabloc
  \def\!tickxsign{0}%
  \def\!tickysign{-}%
  \def\!axissetup{\!axisxsetup}%
  \def\!axislabeltbrl{t}%
  \!nextkeyword}

\def\!axistop{%
  \!axisylevel=\!areatloc
  \def\!tickxsign{0}%
  \def\!tickysign{+}%
  \def\!axissetup{\!axisxsetup}%
  \def\!axislabeltbrl{b}%
  \!nextkeyword}

\def\!axisleft{%
  \!axisxlevel=\!arealloc
  \def\!tickxsign{-}%
  \def\!tickysign{0}%
  \def\!axissetup{\!axisysetup}%
  \def\!axislabeltbrl{r}%
  \!nextkeyword}

\def\!axisright{%
  \!axisxlevel=\!arearloc
  \def\!tickxsign{+}%
  \def\!tickysign{0}%
  \def\!axissetup{\!axisysetup}%
  \def\!axislabeltbrl{l}%
  \!nextkeyword}

\def\!axisshiftedto#1=#2 {%
  \if 0\!tickxsign
    \!axisylevel=\!M{#2}\!yunit
    \advance\!axisylevel -\!yorigin
  \else
    \!axisxlevel=\!M{#2}\!xunit
    \advance\!axisxlevel -\!xorigin
  \fi
  \!nextkeyword}

\def\!axisvisible{%
  \!axisvisibletrue  
  \!nextkeyword}

\def\!axisinvisible{%
  \!axisvisiblefalse
  \!nextkeyword}

\def\!axislabel#1 {%
  \!axisLaBeL={#1}%
  \!placeaxislabeltrue
  \!nextkeyword}

\expandafter\def\csname !axis/\endcsname{%
  \!axissetup 
  \if!placeaxislabel
    \!placeaxislabel
  \fi
  \if +\!tickysign 
    \!dimenA=\!axisylevel
    \advance\!dimenA \!offset 
    \advance\!dimenA -\!areatloc 
    \ifdim \!dimenA>\!plotheadingoffset
      \!plotheadingoffset=\!dimenA 
    \fi
  \fi}

\def\grid #1 #2 {%
  \!countA=#1\advance\!countA 1
  \axis bottom invisible ticks length <\!zpt> andacross quantity {\!countA} /
  \!countA=#2\advance\!countA 1
  \axis left   invisible ticks length <\!zpt> andacross quantity {\!countA} / }

\def\plotheading#1 {%
  \advance\!plotheadingoffset \headingtoplotskip
  \!dimenput {#1} [B] <.5\!xaxislength,\!plotheadingoffset>
    (\!arealloc,\!areatloc)}

\def\!axisxsetup{%
  \!axisxlevel=\!arealloc
  \!axisstart=\!arealloc
  \!axisend=\!arearloc
  \!axisLength=\!xaxislength
  \!!origin=\!xorigin
  \!!unit=\!xunit
  \!xswitchtrue
  \if!axisvisible 
    \!makeaxis
  \fi}

\def\!axisysetup{%
  \!axisylevel=\!areabloc
  \!axisstart=\!areabloc
  \!axisend=\!areatloc
  \!axisLength=\!yaxislength
  \!!origin=\!yorigin
  \!!unit=\!yunit
  \!xswitchfalse
  \if!axisvisible
    \!makeaxis
  \fi}

\def\!makeaxis{%
  \setbox\!boxA=\hbox{
    \beginpicture
      \!setdimenmode
      \setcoordinatesystem point at {\!zpt} {\!zpt}   
      \putrule from {\!zpt} {\!zpt} to
        {\!tickysign\!tickysign\!axisLength} 
        {\!tickxsign\!tickxsign\!axisLength}
    \endpicturesave <\!Xsave,\!Ysave>}%
    \wd\!boxA=\!zpt
    \!placetick\!axisstart}

\def\!placeaxislabel{%
  \advance\!offset \valuestolabelleading
  \if!xswitch
    \!dimenput {\the\!axisLaBeL} [\!axislabeltbrl]
      <.5\!axisLength,\!tickysign\!offset> (\!axisxlevel,\!axisylevel)
    \advance\!offset \!dp  
    \advance\!offset \!ht  
  \else
    \!dimenput {\the\!axisLaBeL} [\!axislabeltbrl]
      <\!tickxsign\!offset,.5\!axisLength> (\!axisxlevel,\!axisylevel)
  \fi
  \!axisLaBeL={}}

\def\arrow <#1> [#2,#3]{%
  \!ifnextchar<{\!arrow{#1}{#2}{#3}}{\!arrow{#1}{#2}{#3}<\!zpt,\!zpt> }}

\def\!arrow#1#2#3<#4,#5> from #6 #7 to #8 #9 {%
  \!xloc=\!M{#8}\!xunit   
  \!yloc=\!M{#9}\!yunit
  \!dxpos=\!xloc  \!dimenA=\!M{#6}\!xunit  \advance \!dxpos -\!dimenA
  \!dypos=\!yloc  \!dimenA=\!M{#7}\!yunit  \advance \!dypos -\!dimenA
  \let\!MAH=\!M
  \!setdimenmode
  \!xshift=#4\relax  \!yshift=#5\relax
  \!reverserotateonly\!xshift\!yshift
  \advance\!xshift\!xloc  \advance\!yshift\!yloc
  \!xS=-\!dxpos  \advance\!xS\!xshift
  \!yS=-\!dypos  \advance\!yS\!yshift
  \!start (\!xS,\!yS)
  \!ljoin (\!xshift,\!yshift)
  \!Pythag\!dxpos\!dypos\!arclength
  \!divide\!dxpos\!arclength\!dxpos  
  \!dxpos=32\!dxpos  \!removept\!dxpos\!!cos
  \!divide\!dypos\!arclength\!dypos  
  \!dypos=32\!dypos  \!removept\!dypos\!!sin
  \!halfhead{#1}{#2}{#3}
  \!halfhead{#1}{-#2}{-#3}
  \let\!M=\!MAH
  \ignorespaces}
  \def\!halfhead#1#2#3{%
    \!dimenC=-#1%
    \divide \!dimenC 2 
    \!dimenD=#2\!dimenC
    \!rotate(\!dimenC,\!dimenD)by(\!!cos,\!!sin)to(\!xM,\!yM)
    \!dimenC=-#1
    \!dimenD=#3\!dimenC
    \!dimenD=.5\!dimenD
    \!rotate(\!dimenC,\!dimenD)by(\!!cos,\!!sin)to(\!xE,\!yE)
    \!start (\!xshift,\!yshift)
    \advance\!xM\!xshift  \advance\!yM\!yshift
    \advance\!xE\!xshift  \advance\!yE\!yshift
    \!qjoin (\!xM,\!yM) (\!xE,\!yE) 
    \ignorespaces}

\def\betweenarrows #1#2 from #3 #4 to #5 #6 {%
  \!xloc=\!M{#3}\!xunit  \!xxloc=\!M{#5}\!xunit%
  \!yloc=\!M{#4}\!yunit  \!yyloc=\!M{#6}\!yunit%
  \!dxpos=\!xxloc  \advance\!dxpos by -\!xloc
  \!dypos=\!yyloc  \advance\!dypos by -\!yloc
  \advance\!xloc .5\!dxpos
  \advance\!yloc .5\!dypos
  \let\!MBA=\!M
  \!setdimenmode
  \ifdim\!dypos=\!zpt
    \ifdim\!dxpos<\!zpt \!dxpos=-\!dxpos \fi
    \put {\!lrarrows{\!dxpos}{#1}}#2{} at {\!xloc} {\!yloc}
  \else
    \ifdim\!dxpos=\!zpt
      \ifdim\!dypos<\!zpt \!dypos=-\!zpt \fi
      \put {\!udarrows{\!dypos}{#1}}#2{} at {\!xloc} {\!yloc}
    \fi
  \fi
  \let\!M=\!MBA
  \ignorespaces}

\def\!lrarrows#1#2{
  {\setbox\!boxA=\hbox{$\mkern-2mu\mathord-\mkern-2mu$}%
   \setbox\!boxB=\hbox{$\leftarrow$}\!dimenE=\ht\!boxB
   \setbox\!boxB=\hbox{}\ht\!boxB=2\!dimenE
   \hbox to #1{$\mathord\leftarrow\mkern-6mu
     \cleaders\copy\!boxA\hfil
     \mkern-6mu\mathord-$%
     \kern.4em $\vcenter{\box\!boxB}$$\vcenter{\hbox{#2}}$\kern.4em
     $\mathord-\mkern-6mu
     \cleaders\copy\!boxA\hfil
     \mkern-6mu\mathord\rightarrow$}}}

\def\!udarrows#1#2{
  {\setbox\!boxB=\hbox{#2}%
   \setbox\!boxA=\hbox to \wd\!boxB{\hss$\vert$\hss}%
   \!dimenE=\ht\!boxA \advance\!dimenE \dp\!boxA \divide\!dimenE 2
   \vbox to #1{\offinterlineskip
      \vskip .05556\!dimenE
      \hbox to \wd\!boxB{\hss$\mkern.4mu\uparrow$\hss}\vskip-\!dimenE
      \cleaders\copy\!boxA\vfil
      \vskip-\!dimenE\copy\!boxA
      \vskip\!dimenE\copy\!boxB\vskip.4em
      \copy\!boxA\vskip-\!dimenE
      \cleaders\copy\!boxA\vfil
      \vskip-\!dimenE \hbox to \wd\!boxB{\hss$\mkern.4mu\downarrow$\hss}
      \vskip .05556\!dimenE}}}

\def\putbar#1breadth <#2> from #3 #4 to #5 #6 {%
  \!xloc=\!M{#3}\!xunit  \!xxloc=\!M{#5}\!xunit%
  \!yloc=\!M{#4}\!yunit  \!yyloc=\!M{#6}\!yunit%
  \!dypos=\!yyloc  \advance\!dypos by -\!yloc
  \!dimenI=#2  
  \ifdim \!dimenI=\!zpt 
    \putrule#1from {#3} {#4} to {#5} {#6} 
  \else 
    \let\!MBar=\!M
    \!setdimenmode 
    \divide\!dimenI 2
    \ifdim \!dypos=\!zpt             
      \advance \!yloc -\!dimenI 
      \advance \!yyloc \!dimenI
    \else
      \advance \!xloc -\!dimenI 
      \advance \!xxloc \!dimenI
    \fi
    \putrectangle#1corners at {\!xloc} {\!yloc} and {\!xxloc} {\!yyloc}
    \let\!M=\!MBar 
  \fi
  \ignorespaces}

\def\setbars#1breadth <#2> baseline at #3 = #4 {%
  \edef\!barshift{#1}%
  \edef\!barbreadth{#2}%
  \edef\!barorientation{#3}%
  \edef\!barbaseline{#4}%
  \def\!bardobaselabel{\!bardoendlabel}%
  \def\!bardoendlabel{\!barfinish}%
  \let\!drawcurve=\!barcurve
  \!setbars}
\def\!setbars{%
  \futurelet\!nextchar\!!setbars}
\def\!!setbars{%
  \if b\!nextchar
    \def\!!!setbars{\!setbarsbget}%
  \else 
    \if e\!nextchar
      \def\!!!setbars{\!setbarseget}%
    \else
      \def\!!!setbars{\relax}%
    \fi
  \fi
  \!!!setbars}
\def\!setbarsbget baselabels (#1) {%
  \def\!barbaselabelorientation{#1}%
  \def\!bardobaselabel{\!!bardobaselabel}%
  \!setbars}
\def\!setbarseget endlabels (#1) {%
  \edef\!barendlabelorientation{#1}%
  \def\!bardoendlabel{\!!bardoendlabel}%
  \!setbars}

\def\!barcurve #1 #2 {%
  \if y\!barorientation
    \def\!basexarg{#1}%
    \def\!baseyarg{\!barbaseline}%
  \else
    \def\!basexarg{\!barbaseline}%
    \def\!baseyarg{#2}%
  \fi
  \expandafter\putbar\!barshift breadth <\!barbreadth> from {\!basexarg}
    {\!baseyarg} to {#1} {#2}
  \def\!endxarg{#1}%
  \def\!endyarg{#2}%
  \!bardobaselabel}

\def\!!bardobaselabel "#1" {%
  \put {#1}\!barbaselabelorientation{} at {\!basexarg} {\!baseyarg}
  \!bardoendlabel}
 
\def\!!bardoendlabel "#1" {%
  \put {#1}\!barendlabelorientation{} at {\!endxarg} {\!endyarg}
  \!barfinish}

\def\!barfinish{%
  \!ifnextchar/{\!finish}{\!barcurve}}

\def\putrectangle{%
  \!ifnextchar<{\!putrectangle}{\!putrectangle<\!zpt,\!zpt> }}
\def\!putrectangle<#1,#2> corners at #3 #4 and #5 #6 {%
  \!xone=\!M{#3}\!xunit  \!xtwo=\!M{#5}\!xunit%
  \!yone=\!M{#4}\!yunit  \!ytwo=\!M{#6}\!yunit%
  \ifdim \!xtwo<\!xone
    \!dimenI=\!xone  \!xone=\!xtwo  \!xtwo=\!dimenI
  \fi
  \ifdim \!ytwo<\!yone
    \!dimenI=\!yone  \!yone=\!ytwo  \!ytwo=\!dimenI
  \fi
  \!dimenI=#1\relax  \advance\!xone\!dimenI  \advance\!xtwo\!dimenI
  \!dimenI=#2\relax  \advance\!yone\!dimenI  \advance\!ytwo\!dimenI
  \let\!MRect=\!M
  \!setdimenmode
  \!shaderectangle
  \!dimenI=.5\linethickness
  \advance \!xone  -\!dimenI
  \advance \!xtwo   \!dimenI
  \putrule from {\!xone} {\!yone} to {\!xtwo} {\!yone} 
  \putrule from {\!xone} {\!ytwo} to {\!xtwo} {\!ytwo} 
  \advance \!xone   \!dimenI
  \advance \!xtwo  -\!dimenI%
  \advance \!yone  -\!dimenI
  \advance \!ytwo   \!dimenI
  \putrule from {\!xone} {\!yone} to {\!xone} {\!ytwo} 
  \putrule from {\!xtwo} {\!yone} to {\!xtwo} {\!ytwo} 
  \let\!M=\!MRect
  \ignorespaces}
 
\def\shaderectanglesoff{%
  \def\!shaderectangle{}%
  \ignorespaces}

\shaderectanglesoff
 
\def\!!shaderectangle{%
  \!dimenA=\!xtwo  \advance \!dimenA -\!xone
  \!dimenB=\!ytwo  \advance \!dimenB -\!yone
  \ifdim \!dimenA<\!dimenB
    \!startvshade (\!xone,\!yone,\!ytwo)
    \!lshade      (\!xtwo,\!yone,\!ytwo)
  \else
    \!starthshade (\!yone,\!xone,\!xtwo)
    \!lshade      (\!ytwo,\!xone,\!xtwo)
  \fi
  \ignorespaces}
  
\def\frame{%
  \!ifnextchar<{\!frame}{\!frame<\!zpt> }}
\long\def\!frame<#1> #2{%
  \beginpicture
    \setcoordinatesystem units <1pt,1pt> point at 0 0 
    \put {#2} [Bl] at 0 0 
    \!dimenA=#1\relax
    \!dimenB=\!wd \advance \!dimenB \!dimenA
    \!dimenC=\!ht \advance \!dimenC \!dimenA
    \!dimenD=\!dp \advance \!dimenD \!dimenA
    \let\!MFr=\!M
    \!setdimenmode
    \putrectangle corners at {-\!dimenA} {-\!dimenD} and {\!dimenB} {\!dimenC}
    \!setcoordmode
    \let\!M=\!MFr
  \endpicture
  \ignorespaces}
 
\def\rectangle <#1> <#2> {%
  \setbox0=\hbox{}\wd0=#1\ht0=#2\frame {\box0}}

\def\!plotfromfile"#1"{%
  \expandafter\!drawcurve \input #1 /}

\def\setquadratic{%
  \let\!drawcurve=\!qcurve
  \let\!!Shade=\!!qShade
  \let\!!!Shade=\!!!qShade}

\def\setlinear{%
  \let\!drawcurve=\!lcurve
  \let\!!Shade=\!!lShade
  \let\!!!Shade=\!!!lShade}

\def\sethistograms{%
  \let\!drawcurve=\!hcurve}

\def\!qcurve #1 #2 {%
  \!start (#1,#2)
  \!Qjoin}
\def\!Qjoin#1 #2 #3 #4 {%
  \!qjoin (#1,#2) (#3,#4)             
  \!ifnextchar/{\!finish}{\!Qjoin}}

\def\!lcurve #1 #2 {%
  \!start (#1,#2)
  \!Ljoin}
\def\!Ljoin#1 #2 {%
  \!ljoin (#1,#2)                    
  \!ifnextchar/{\!finish}{\!Ljoin}}

\def\!finish/{\ignorespaces}

\def\!hcurve #1 #2 {%
  \edef\!hxS{#1}%
  \edef\!hyS{#2}%
  \!hjoin}
\def\!hjoin#1 #2 {%
  \putrectangle corners at {\!hxS} {\!hyS} and {#1} {#2}
  \edef\!hxS{#1}%
  \!ifnextchar/{\!finish}{\!hjoin}}

\def\vshade #1 #2 #3 {%
  \!startvshade (#1,#2,#3)
  \!Shadewhat}

\def\hshade #1 #2 #3 {%
  \!starthshade (#1,#2,#3)
  \!Shadewhat}

\def\!Shadewhat{%
  \futurelet\!nextchar\!Shade}
\def\!Shade{%
  \if <\!nextchar
    \def\!nextShade{\!!Shade}%
  \else
    \if /\!nextchar
      \def\!nextShade{\!finish}%
    \else
      \def\!nextShade{\!!!Shade}%
    \fi
  \fi
  \!nextShade}
\def\!!lShade<#1> #2 #3 #4 {%
  \!lshade <#1> (#2,#3,#4)                 
  \!Shadewhat}
\def\!!!lShade#1 #2 #3 {%
  \!lshade (#1,#2,#3)
  \!Shadewhat} 
\def\!!qShade<#1> #2 #3 #4 #5 #6 #7 {%
  \!qshade <#1> (#2,#3,#4) (#5,#6,#7)      
  \!Shadewhat}
\def\!!!qShade#1 #2 #3 #4 #5 #6 {%
  \!qshade (#1,#2,#3) (#4,#5,#6)
  \!Shadewhat} 

\setlinear

\def\setdashpattern <#1>{%
  \def\!Flist{}\def\!Blist{}\def\!UDlist{}%
  \!countA=0
  \!ecfor\!item:=#1\do{%
    \!dimenA=\!item\relax
    \expandafter\!rightappend\the\!dimenA\withCS{\\}\to\!UDlist%
    \advance\!countA  1
    \ifodd\!countA
      \expandafter\!rightappend\the\!dimenA\withCS{\!Rule}\to\!Flist%
      \expandafter\!leftappend\the\!dimenA\withCS{\!Rule}\to\!Blist%
    \else 
      \expandafter\!rightappend\the\!dimenA\withCS{\!Skip}\to\!Flist%
      \expandafter\!leftappend\the\!dimenA\withCS{\!Skip}\to\!Blist%
    \fi}%
  \!leaderlength=\!zpt
  \def\!Rule##1{\advance\!leaderlength  ##1}%
  \def\!Skip##1{\advance\!leaderlength  ##1}%
  \!Flist%
  \ifdim\!leaderlength>\!zpt 
  \else
    \def\!Flist{\!Skip{24in}}\def\!Blist{\!Skip{24in}}\ignorespaces
    \def\!UDlist{\\{\!zpt}\\{24in}}\ignorespaces
    \!leaderlength=24in
  \fi
  \!dashingon}

\def\!dashingon{%
  \def\!advancedashing{\!!advancedashing}%
  \def\!drawlinearsegment{\!lineardashed}%
  \def\!puthline{\!putdashedhline}%
  \def\!putvline{\!putdashedvline}%
  \ignorespaces}%
\def\!dashingoff{%
  \def\!advancedashing{\relax}%
  \def\!drawlinearsegment{\!linearsolid}%
  \def\!puthline{\!putsolidhline}%
  \def\!putvline{\!putsolidvline}%
  \ignorespaces}

\def\setdots{%
  \!ifnextchar<{\!setdots}{\!setdots<5pt>}}
\def\!setdots<#1>{%
  \!dimenB=#1\advance\!dimenB -\plotsymbolspacing
  \ifdim\!dimenB<\!zpt
    \!dimenB=\!zpt
  \fi
\setdashpattern <\plotsymbolspacing,\!dimenB>}
 
\def\setdotsnear <#1> for <#2>{%
  \!dimenB=#2\relax  \advance\!dimenB -.05pt  
  \!dimenC=#1\relax  \!countA=\!dimenC 
  \!dimenD=\!dimenB  \advance\!dimenD .5\!dimenC  \!countB=\!dimenD
  \divide \!countB  \!countA
  \ifnum 1>\!countB 
    \!countB=1
  \fi
  \divide\!dimenB  \!countB
  \setdots <\!dimenB>}
 
\def\setdashes{%
  \!ifnextchar<{\!setdashes}{\!setdashes<5pt>}}
\def\!setdashes<#1>{\setdashpattern <#1,#1>}
 
\def\setdashesnear <#1> for <#2>{%
  \!dimenB=#2\relax  
  \!dimenC=#1\relax  \!countA=\!dimenC 
  \!dimenD=\!dimenB  \advance\!dimenD .5\!dimenC  \!countB=\!dimenD
  \divide \!countB  \!countA
  \ifodd \!countB 
  \else 
    \advance \!countB  1
  \fi
  \divide\!dimenB  \!countB
  \setdashes <\!dimenB>}
 
\def\setsolid{%
  \def\!Flist{\!Rule{24in}}\def\!Blist{\!Rule{24in}}%
  \def\!UDlist{\\{24in}\\{\!zpt}}%
  \!dashingoff}  
\setsolid

\def\!divide#1#2#3{%
  \!dimenB=#1
  \!dimenC=#2
  \!dimenD=\!dimenB
  \divide \!dimenD \!dimenC
  \!dimenA=\!dimenD
  \multiply\!dimenD \!dimenC
  \advance\!dimenB -\!dimenD
  \!dimenD=\!dimenC
    \ifdim\!dimenD<\!zpt \!dimenD=-\!dimenD 
  \fi
  \ifdim\!dimenD<64pt
    \!divstep[\!tfs]\!divstep[\!tfs]%
  \else 
    \!!divide
  \fi
  #3=\!dimenA\ignorespaces}

\def\!!divide{%
  \ifdim\!dimenD<256pt
    \!divstep[64]\!divstep[32]\!divstep[32]%
  \else 
    \!divstep[8]\!divstep[8]\!divstep[8]\!divstep[8]\!divstep[8]%
    \!dimenA=2\!dimenA
  \fi}

\def\!divstep[#1]{
  \!dimenB=#1\!dimenB
  \!dimenD=\!dimenB
    \divide \!dimenD by \!dimenC
  \!dimenA=#1\!dimenA
    \advance\!dimenA by \!dimenD%
  \multiply\!dimenD by \!dimenC
    \advance\!dimenB by -\!dimenD}
 
\def\Divide <#1> by <#2> forming <#3> {%
  \!divide{#1}{#2}{#3}}

\def\ellipticalarc axes ratio #1:#2 #3 degrees from #4 #5 center at #6 #7 {%
  \!angle=#3pt\relax
  \ifdim\!angle>\!zpt 
    \def\!sign{}
  \else 
    \def\!sign{-}\!angle=-\!angle
  \fi
  \!xxloc=\!M{#6}\!xunit
  \!yyloc=\!M{#7}\!yunit     
  \!xxS=\!M{#4}\!xunit
  \!yyS=\!M{#5}\!yunit
  \advance\!xxS -\!xxloc
  \advance\!yyS -\!yyloc
  \!divide\!xxS{#1pt}\!xxS 
  \!divide\!yyS{#2pt}\!yyS 
  \let\!MC=\!M
  \!setdimenmode
  \!xS=#1\!xxS  \advance\!xS\!xxloc
  \!yS=#2\!yyS  \advance\!yS\!yyloc
  \!start (\!xS,\!yS)%
  \!loop\ifdim\!angle>14.9999pt
    \!rotate(\!xxS,\!yyS)by(\!cos,\!sign\!sin)to(\!xxM,\!yyM) 
    \!rotate(\!xxM,\!yyM)by(\!cos,\!sign\!sin)to(\!xxE,\!yyE)
    \!xM=#1\!xxM  \advance\!xM\!xxloc  \!yM=#2\!yyM  \advance\!yM\!yyloc
    \!xE=#1\!xxE  \advance\!xE\!xxloc  \!yE=#2\!yyE  \advance\!yE\!yyloc
    \!qjoin (\!xM,\!yM) (\!xE,\!yE)
    \!xxS=\!xxE  \!yyS=\!yyE 
    \advance \!angle -15pt
  \repeat
  \ifdim\!angle>\!zpt
    \!angle=100.53096\!angle
    \divide \!angle 360 
    \!sinandcos\!angle\!!sin\!!cos
    \!rotate(\!xxS,\!yyS)by(\!!cos,\!sign\!!sin)to(\!xxM,\!yyM) 
    \!rotate(\!xxM,\!yyM)by(\!!cos,\!sign\!!sin)to(\!xxE,\!yyE)
    \!xM=#1\!xxM  \advance\!xM\!xxloc  \!yM=#2\!yyM  \advance\!yM\!yyloc
    \!xE=#1\!xxE  \advance\!xE\!xxloc  \!yE=#2\!yyE  \advance\!yE\!yyloc
    \!qjoin (\!xM,\!yM) (\!xE,\!yE)
  \fi
  \let\!M=\!MC
  \ignorespaces}

\def\!rotate(#1,#2)by(#3,#4)to(#5,#6){%
  \!dimenA=#3#1\advance \!dimenA -#4#2
  \!dimenB=#3#2\advance \!dimenB  #4#1
  \divide \!dimenA 32  \divide \!dimenB 32 
  #5=\!dimenA  #6=\!dimenB
  \ignorespaces}
\def\!sin{4.17684}
\def\!cos{31.72624}

\def\!sinandcos#1#2#3{%
 \!dimenD=#1
 \!dimenA=\!dimenD
 \!dimenB=32pt
 \!removept\!dimenD\!value
 \!dimenC=\!dimenD
 \!dimenC=\!value\!dimenC \divide\!dimenC by 64 
 \advance\!dimenB by -\!dimenC
 \!dimenC=\!value\!dimenC \divide\!dimenC by 96 
 \advance\!dimenA by -\!dimenC
 \!dimenC=\!value\!dimenC \divide\!dimenC by 128 
 \advance\!dimenB by \!dimenC%
 \!removept\!dimenA#2
 \!removept\!dimenB#3
 \ignorespaces}

\def\putrule#1from #2 #3 to #4 #5 {%
  \!xloc=\!M{#2}\!xunit  \!xxloc=\!M{#4}\!xunit%
  \!yloc=\!M{#3}\!yunit  \!yyloc=\!M{#5}\!yunit%
  \!dxpos=\!xxloc  \advance\!dxpos by -\!xloc
  \!dypos=\!yyloc  \advance\!dypos by -\!yloc
  \ifdim\!dypos=\!zpt
    \def\!!Line{\!puthline{#1}}\ignorespaces
  \else
    \ifdim\!dxpos=\!zpt
      \def\!!Line{\!putvline{#1}}\ignorespaces
    \else 
       \def\!!Line{}
    \fi
  \fi
  \let\!ML=\!M
  \!setdimenmode
  \!!Line%
  \let\!M=\!ML
  \ignorespaces}

\def\!putsolidhline#1{%
  \ifdim\!dxpos>\!zpt 
    \put{\!hline\!dxpos}#1[l] at {\!xloc} {\!yloc}
  \else 
    \put{\!hline{-\!dxpos}}#1[l] at {\!xxloc} {\!yyloc}
  \fi
  \ignorespaces}
 
\def\!putsolidvline#1{%
  \ifdim\!dypos>\!zpt 
    \put{\!vline\!dypos}#1[b] at {\!xloc} {\!yloc}
  \else 
    \put{\!vline{-\!dypos}}#1[b] at {\!xxloc} {\!yyloc}
  \fi
  \ignorespaces}
 
\def\!hline#1{\hbox to #1{\leaders \hrule height\linethickness\hfill}}
\def\!vline#1{\vbox to #1{\leaders \vrule width\linethickness\vfill}}

\def\!putdashedhline#1{%
  \ifdim\!dxpos>\!zpt 
    \!DLsetup\!Flist\!dxpos
    \put{\hbox to \!totalleaderlength{\!hleaders}\!hpartialpattern\!Rtrunc}
      #1[l] at {\!xloc} {\!yloc} 
  \else 
    \!DLsetup\!Blist{-\!dxpos}
    \put{\!hpartialpattern\!Ltrunc\hbox to \!totalleaderlength{\!hleaders}}
      #1[r] at {\!xloc} {\!yloc} 
  \fi
  \ignorespaces}
 
\def\!putdashedvline#1{%
  \!dypos=-\!dypos
  \ifdim\!dypos>\!zpt 
    \!DLsetup\!Flist\!dypos 
    \put{\vbox{\vbox to \!totalleaderlength{\!vleaders}
      \!vpartialpattern\!Rtrunc}}#1[t] at {\!xloc} {\!yloc} 
  \else 
    \!DLsetup\!Blist{-\!dypos}
    \put{\vbox{\!vpartialpattern\!Ltrunc
      \vbox to \!totalleaderlength{\!vleaders}}}#1[b] at {\!xloc} {\!yloc} 
  \fi
  \ignorespaces}

\def\!DLsetup#1#2{
  \let\!RSlist=#1
  \!countB=#2
  \!countA=\!leaderlength
  \divide\!countB by \!countA
  \!totalleaderlength=\!countB\!leaderlength
  \!Rresiduallength=#2%
  \advance \!Rresiduallength by -\!totalleaderlength
  \!Lresiduallength=\!leaderlength
  \advance \!Lresiduallength by -\!Rresiduallength
  \ignorespaces}
 
\def\!hleaders{%
  \def\!Rule##1{\vrule height\linethickness width##1}%
  \def\!Skip##1{\hskip##1}%
  \leaders\hbox{\!RSlist}\hfill}
 
\def\!hpartialpattern#1{%
  \!dimenA=\!zpt \!dimenB=\!zpt 
  \def\!Rule##1{#1{##1}\vrule height\linethickness width\!dimenD}%
  \def\!Skip##1{#1{##1}\hskip\!dimenD}%
  \!RSlist}
 
\def\!vleaders{%
  \def\!Rule##1{\hrule width\linethickness height##1}%
  \def\!Skip##1{\vskip##1}%
  \leaders\vbox{\!RSlist}\vfill}
 
\def\!vpartialpattern#1{%
  \!dimenA=\!zpt \!dimenB=\!zpt 
  \def\!Rule##1{#1{##1}\hrule width\linethickness height\!dimenD}%
  \def\!Skip##1{#1{##1}\vskip\!dimenD}%
  \!RSlist}
 
\def\!Rtrunc#1{\!trunc{#1}>\!Rresiduallength}
\def\!Ltrunc#1{\!trunc{#1}<\!Lresiduallength}
 
\def\!trunc#1#2#3{%
  \!dimenA=\!dimenB         
  \advance\!dimenB by #1%
  \!dimenD=\!dimenB  \ifdim\!dimenD#2#3\!dimenD=#3\fi
  \!dimenC=\!dimenA  \ifdim\!dimenC#2#3\!dimenC=#3\fi
  \advance \!dimenD by -\!dimenC}

\def\!start (#1,#2){%
  \!plotxorigin=\!xorigin  \advance \!plotxorigin by \!plotsymbolxshift
  \!plotyorigin=\!yorigin  \advance \!plotyorigin by \!plotsymbolyshift
  \!xS=\!M{#1}\!xunit \!yS=\!M{#2}\!yunit
  \!rotateaboutpivot\!xS\!yS
  \!copylist\!UDlist\to\!!UDlist
  \!getnextvalueof\!downlength\from\!!UDlist
  \!distacross=\!zpt
  \!intervalno=0 
  \global\totalarclength=\!zpt
  \ignorespaces}

\def\!ljoin (#1,#2){%
  \advance\!intervalno by 1
  \!xE=\!M{#1}\!xunit \!yE=\!M{#2}\!yunit
  \!rotateaboutpivot\!xE\!yE
  \!xdiff=\!xE \advance \!xdiff by -\!xS
  \!ydiff=\!yE \advance \!ydiff by -\!yS
  \!Pythag\!xdiff\!ydiff\!arclength
  \global\advance \totalarclength by \!arclength%
  \!drawlinearsegment
  \!xS=\!xE \!yS=\!yE
  \ignorespaces}

\def\!linearsolid{%
  \!npoints=\!arclength
  \!countA=\plotsymbolspacing
  \divide\!npoints by \!countA
  \ifnum \!npoints<1 
    \!npoints=1 
  \fi
  \divide\!xdiff by \!npoints
  \divide\!ydiff by \!npoints
  \!xpos=\!xS \!ypos=\!yS
  \loop\ifnum\!npoints>-1
    \!plotifinbounds
    \advance \!xpos by \!xdiff
    \advance \!ypos by \!ydiff
    \advance \!npoints by -1
  \repeat
  \ignorespaces}

\def\!lineardashed{%
  \ifdim\!distacross>\!arclength
    \advance \!distacross by -\!arclength  
  \else
    \loop\ifdim\!distacross<\!arclength
      \!divide\!distacross\!arclength\!dimenA
      \!removept\!dimenA\!t
      \!xpos=\!t\!xdiff \advance \!xpos by \!xS
      \!ypos=\!t\!ydiff \advance \!ypos by \!yS
      \!plotifinbounds
      \advance\!distacross by \plotsymbolspacing
      \!advancedashing
    \repeat  
    \advance \!distacross by -\!arclength
  \fi
  \ignorespaces}

\def\!!advancedashing{%
  \advance\!downlength by -\plotsymbolspacing
  \ifdim \!downlength>\!zpt
  \else
    \advance\!distacross by \!downlength
    \!getnextvalueof\!uplength\from\!!UDlist
    \advance\!distacross by \!uplength
    \!getnextvalueof\!downlength\from\!!UDlist
  \fi}

\def\inboundscheckoff{%
  \def\!plotifinbounds{\!plot(\!xpos,\!ypos)}%
  \def\!initinboundscheck{\relax}\ignorespaces}
 
\inboundscheckoff
 
\def\!!plotifinbounds{%
  \ifdim \!xpos<\!checkleft
  \else
    \ifdim \!xpos>\!checkright
    \else
      \ifdim \!ypos<\!checkbot
      \else
         \ifdim \!ypos>\!checktop
         \else
           \!plot(\!xpos,\!ypos)
         \fi 
      \fi
    \fi
  \fi}

\def\!!initinboundscheck{%
  \!checkleft=\!arealloc     \advance\!checkleft by \!xorigin
  \!checkright=\!arearloc    \advance\!checkright by \!xorigin
  \!checkbot=\!areabloc      \advance\!checkbot by \!yorigin
  \!checktop=\!areatloc      \advance\!checktop by \!yorigin}

\def\!logten#1#2{%
  \expandafter\!!logten#1\!nil
  \!removept\!dimenF#2%
  \ignorespaces}

\def\!!logten#1#2\!nil{%
  \if -#1%
    \!dimenF=\!zpt
    \def\!next{\ignorespaces}%
  \else
    \if +#1%
      \def\!next{\!!logten#2\!nil}%
    \else
      \if .#1%
        \def\!next{\!!logten0.#2\!nil}%
      \else
        \def\!next{\!!!logten#1#2..\!nil}%
      \fi
    \fi
  \fi
  \!next}

\def\!!!logten#1#2.#3.#4\!nil{%
  \!dimenF=1pt 
  \if 0#1%
    \!!logshift#3pt 
  \else 
    \!logshift#2/
    \!dimenE=#1.#2#3pt 
  \fi 
  \ifdim \!dimenE<\!rootten
    \multiply \!dimenE 10 
    \advance  \!dimenF -1pt
  \fi
  \!dimenG=\!dimenE
    \advance\!dimenG 10pt
  \advance\!dimenE -10pt 
  \multiply\!dimenE 10 
  \!divide\!dimenE\!dimenG\!dimenE
  \!removept\!dimenE\!t
  \!dimenG=\!t\!dimenE
  \!removept\!dimenG\!tt
  \!dimenH=\!tt\!tenAe
    \divide\!dimenH 100
  \advance\!dimenH \!tenAc
  \!dimenH=\!tt\!dimenH
    \divide\!dimenH 100   
  \advance\!dimenH \!tenAa
  \!dimenH=\!t\!dimenH
    \divide\!dimenH 100 
  \advance\!dimenF \!dimenH}

\def\!logshift#1{%
  \if #1/%
    \def\!next{\ignorespaces}%
  \else
    \advance\!dimenF 1pt 
    \def\!next{\!logshift}%
  \fi 
  \!next}
 
 \def\!!logshift#1{%
   \advance\!dimenF -1pt
   \if 0#1%
     \def\!next{\!!logshift}%
   \else
     \if p#1%
       \!dimenF=1pt
       \def\!next{\!dimenE=1p}%
     \else
       \def\!next{\!dimenE=#1.}%
     \fi
   \fi
   \!next}

\def\beginpicture{%
  \setbox\!picbox=\hbox\bgroup%
  \!xleft=\maxdimen  
  \!xright=-\maxdimen
  \!ybot=\maxdimen
  \!ytop=-\maxdimen}
 
\def\endpicture{%
  \ifdim\!xleft=\maxdimen
    \!xleft=\!zpt \!xright=\!zpt \!ybot=\!zpt \!ytop=\!zpt 
  \fi
  \global\!Xleft=\!xleft \global\!Xright=\!xright
  \global\!Ybot=\!ybot \global\!Ytop=\!ytop
  \egroup%
  \ht\!picbox=\!Ytop  \dp\!picbox=-\!Ybot
  \ifdim\!Ybot>\!zpt
  \else 
    \ifdim\!Ytop<\!zpt
      \!Ybot=\!Ytop
    \else
      \!Ybot=\!zpt
    \fi
  \fi
  \hbox{\kern-\!Xleft\lower\!Ybot\box\!picbox\kern\!Xright}}
 
\def\endpicturesave <#1,#2>{%
  \endpicture \global #1=\!Xleft \global #2=\!Ybot \ignorespaces}

\def\setcoordinatesystem{%
  \!ifnextchar{u}{\!getlengths }
    {\!getlengths units <\!xunit,\!yunit>}}
\def\!getlengths units <#1,#2>{%
  \!xunit=#1\relax
  \!yunit=#2\relax
  \!ifcoordmode 
    \let\!SCnext=\!SCccheckforRP
  \else
    \let\!SCnext=\!SCdcheckforRP
  \fi
  \!SCnext}
\def\!SCccheckforRP{%
  \!ifnextchar{p}{\!cgetreference }
    {\!cgetreference point at {\!xref} {\!yref} }}
\def\!cgetreference point at #1 #2 {%
  \edef\!xref{#1}\edef\!yref{#2}%
  \!xorigin=\!xref\!xunit  \!yorigin=\!yref\!yunit  
  \!initinboundscheck 
  \ignorespaces}
\def\!SCdcheckforRP{%
  \!ifnextchar{p}{\!dgetreference}%
    {\ignorespaces}}
\def\!dgetreference point at #1 #2 {%
  \!xorigin=#1\relax  \!yorigin=#2\relax
  \ignorespaces}

\long\def\put#1#2 at #3 #4 {%
  \!setputobject{#1}{#2}%
  \!xpos=\!M{#3}\!xunit  \!ypos=\!M{#4}\!yunit  
  \!rotateaboutpivot\!xpos\!ypos%
  \advance\!xpos -\!xorigin  \advance\!xpos -\!xshift
  \advance\!ypos -\!yorigin  \advance\!ypos -\!yshift
  \kern\!xpos\raise\!ypos\box\!putobject\kern-\!xpos%
  \!doaccounting\ignorespaces}
 
\long\def\multiput #1#2 at {%
  \!setputobject{#1}{#2}%
  \!ifnextchar"{\!putfromfile}{\!multiput}}
\def\!putfromfile"#1"{%
  \expandafter\!multiput \input #1 /}
\def\!multiput{%
  \futurelet\!nextchar\!!multiput}
\def\!!multiput{%
  \if *\!nextchar
    \def\!nextput{\!alsoby}%
  \else
    \if /\!nextchar
      \def\!nextput{\!finishmultiput}%
    \else
      \def\!nextput{\!alsoat}%
    \fi
  \fi
  \!nextput}
\def\!finishmultiput/{%
  \setbox\!putobject=\hbox{}%
  \ignorespaces}
 
\def\!alsoat#1 #2 {%
  \!xpos=\!M{#1}\!xunit  \!ypos=\!M{#2}\!yunit  
  \!rotateaboutpivot\!xpos\!ypos%
  \advance\!xpos -\!xorigin  \advance\!xpos -\!xshift
  \advance\!ypos -\!yorigin  \advance\!ypos -\!yshift
  \kern\!xpos\raise\!ypos\copy\!putobject\kern-\!xpos%
  \!doaccounting
  \!multiput}
 
\def\!alsoby*#1 #2 #3 {%
  \!dxpos=\!M{#2}\!xunit \!dypos=\!M{#3}\!yunit 
  \!rotateonly\!dxpos\!dypos
  \!ntemp=#1%
  \!!loop\ifnum\!ntemp>0
    \advance\!xpos by \!dxpos  \advance\!ypos by \!dypos
    \kern\!xpos\raise\!ypos\copy\!putobject\kern-\!xpos%
    \advance\!ntemp by -1
  \repeat
  \!doaccounting 
  \!multiput}
 
\def\accountingon{\def\!doaccounting{\!!doaccounting}\ignorespaces}

\accountingon
\def\!!doaccounting{%
  \!xtemp=\!xpos  
  \!ytemp=\!ypos
  \ifdim\!xtemp<\!xleft 
     \!xleft=\!xtemp 
  \fi
  \advance\!xtemp by  \!wd 
  \ifdim\!xright<\!xtemp 
    \!xright=\!xtemp
  \fi
  \advance\!ytemp by -\!dp
  \ifdim\!ytemp<\!ybot  
    \!ybot=\!ytemp
  \fi
  \advance\!ytemp by  \!dp
  \advance\!ytemp by  \!ht 
  \ifdim\!ytemp>\!ytop  
    \!ytop=\!ytemp  
  \fi}
 
\long\def\!setputobject#1#2{%
  \setbox\!putobject=\hbox{#1}%
  \!ht=\ht\!putobject  \!dp=\dp\!putobject  \!wd=\wd\!putobject
  \wd\!putobject=\!zpt
  \!xshift=.5\!wd   \!yshift=.5\!ht   \advance\!yshift by -.5\!dp
  \edef\!putorientation{#2}%
  \expandafter\!SPOreadA\!putorientation[]\!nil%
  \expandafter\!SPOreadB\!putorientation<\!zpt,\!zpt>\!nil\ignorespaces}
 
\def\!SPOreadA#1[#2]#3\!nil{\!etfor\!orientation:=#2\do\!SPOreviseshift}
 
\def\!SPOreadB#1<#2,#3>#4\!nil{\advance\!xshift by -#2\advance\!yshift by -#3}
 
\def\!SPOreviseshift{%
  \if l\!orientation 
    \!xshift=\!zpt
  \else 
    \if r\!orientation 
      \!xshift=\!wd
    \else 
      \if b\!orientation
        \!yshift=-\!dp
      \else 
        \if B\!orientation 
          \!yshift=\!zpt
        \else 
          \if t\!orientation 
            \!yshift=\!ht
          \fi 
        \fi
      \fi
    \fi
  \fi}

\long\def\!dimenput#1#2(#3,#4){%
  \!setputobject{#1}{#2}%
  \!xpos=#3\advance\!xpos by -\!xshift
  \!ypos=#4\advance\!ypos by -\!yshift
  \kern\!xpos\raise\!ypos\box\!putobject\kern-\!xpos%
  \!doaccounting\ignorespaces}

\def\!setdimenmode{%
  \let\!M=\!M!!\ignorespaces}
\def\!setcoordmode{%
  \let\!M=\!M!\ignorespaces}
\def\!ifcoordmode{%
  \ifx \!M \!M!}
\def\!ifdimenmode{%
  \ifx \!M \!M!!}
\def\!M!#1#2{#1#2} 
\def\!M!!#1#2{#1}
\!setcoordmode
\let\setdimensionmode=\!setdimenmode
\let\setcoordinatemode=\!setcoordmode

\def\!stack[#1]{%
  \let\!lglue=\hfill \let\!rglue=\hfill
  \expandafter\let\csname !#1glue\endcsname=\relax
  \!ifnextchar<{\!!stack}{\!!stack<\stackleading>}}
\def\!!stack<#1>#2{%
  \vbox{\def\!valueslist{}\!ecfor\!value:=#2\do{%
    \expandafter\!rightappend\!value\withCS{\\}\to\!valueslist}%
    \!lop\!valueslist\to\!value
    \let\\=\cr\lineskiplimit=\maxdimen\lineskip=#1%
    \baselineskip=-1000pt\halign{\!lglue##\!rglue\cr \!value\!valueslist\cr}}%
  \ignorespaces}

\def\!lines[#1]#2{%
  \let\!lglue=\hfill \let\!rglue=\hfill
  \expandafter\let\csname !#1glue\endcsname=\relax
  \vbox{\halign{\!lglue##\!rglue\cr #2\crcr}}%
  \ignorespaces}

\def\!Lines[#1]#2{%
  \let\!lglue=\hfill \let\!rglue=\hfill
  \expandafter\let\csname !#1glue\endcsname=\relax
  \vtop{\halign{\!lglue##\!rglue\cr #2\crcr}}%
  \ignorespaces}

\def\setplotsymbol(#1#2){%
  \!setputobject{#1}{#2}
  \setbox\!plotsymbol=\box\!putobject%
  \!plotsymbolxshift=\!xshift 
  \!plotsymbolyshift=\!yshift 
  \ignorespaces}
 
\setplotsymbol({\fiverm .})

\def\!!plot(#1,#2){%
  \!dimenA=-\!plotxorigin \advance \!dimenA by #1
  \!dimenB=-\!plotyorigin \advance \!dimenB by #2
  \kern\!dimenA\raise\!dimenB\copy\!plotsymbol\kern-\!dimenA%
  \ignorespaces}
 
\def\!!!plot(#1,#2){%
  \!dimenA=-\!plotxorigin \advance \!dimenA by #1
  \!dimenB=-\!plotyorigin \advance \!dimenB by #2
  \kern\!dimenA\raise\!dimenB\copy\!plotsymbol\kern-\!dimenA%
  \!countE=\!dimenA
  \!countF=\!dimenB
  \immediate\write\!replotfile{\the\!countE,\the\!countF.}%
  \ignorespaces}

\def\savelinesandcurves on "#1" {%
  \immediate\closeout\!replotfile
  \immediate\openout\!replotfile=#1%
  \let\!plot=\!!!plot}

\def\dontsavelinesandcurves {%
  \let\!plot=\!!plot}
\dontsavelinesandcurves

{\catcode`\%=11\xdef\!Commentsignal{
\def\writesavefile#1 {%
  \immediate\write\!replotfile{\!Commentsignal #1}%
  \ignorespaces}

\def\replot"#1" {%
  \expandafter\!replot\input #1 /}
\def\!replot#1,#2. {%
  \!dimenA=#1sp
  \kern\!dimenA\raise#2sp\copy\!plotsymbol\kern-\!dimenA
  \futurelet\!nextchar\!!replot}
\def\!!replot{%
  \if /\!nextchar 
    \def\!next{\!finish}%
  \else
    \def\!next{\!replot}%
  \fi
  \!next}

\def\!Pythag#1#2#3{%
  \!dimenE=#1\relax                                     
  \ifdim\!dimenE<\!zpt 
    \!dimenE=-\!dimenE 
  \fi
  \!dimenF=#2\relax
  \ifdim\!dimenF<\!zpt 
    \!dimenF=-\!dimenF 
  \fi
  \advance \!dimenF by \!dimenE
  \ifdim\!dimenF=\!zpt 
    \!dimenG=\!zpt
  \else 
    \!divide{8\!dimenE}\!dimenF\!dimenE
    \advance\!dimenE by -4pt
      \!dimenE=2\!dimenE
    \!removept\!dimenE\!!t
    \!dimenE=\!!t\!dimenE
    \advance\!dimenE by 64pt
    \divide \!dimenE by 2
    \!dimenH=7pt
    \!!Pythag\!!Pythag\!!Pythag
    \!removept\!dimenH\!!t
    \!dimenG=\!!t\!dimenF
    \divide\!dimenG by 8
  \fi
  #3=\!dimenG
  \ignorespaces}

\def\!!Pythag{
  \!divide\!dimenE\!dimenH\!dimenI
  \advance\!dimenH by \!dimenI
    \divide\!dimenH by 2}

\def\placehypotenuse for <#1> and <#2> in <#3> {%
  \!Pythag{#1}{#2}{#3}}

\def\!qjoin (#1,#2) (#3,#4){%
  \advance\!intervalno by 1
  \!ifcoordmode
    \edef\!xmidpt{#1}\edef\!ymidpt{#2}%
  \else
    \!dimenA=#1\relax \edef\!xmidpt{\the\!dimenA}%
    \!dimenA=#2\relax \edef\!xmidpt{\the\!dimenA}%
  \fi
  \!xM=\!M{#1}\!xunit  \!yM=\!M{#2}\!yunit   \!rotateaboutpivot\!xM\!yM
  \!xE=\!M{#3}\!xunit  \!yE=\!M{#4}\!yunit   \!rotateaboutpivot\!xE\!yE
  \!dimenA=\!xM  \advance \!dimenA by -\!xS
  \!dimenB=\!xE  \advance \!dimenB by -\!xM
  \!xB=3\!dimenA \advance \!xB by -\!dimenB
  \!xC=2\!dimenB \advance \!xC by -2\!dimenA
  \!dimenA=\!yM  \advance \!dimenA by -\!yS%
  \!dimenB=\!yE  \advance \!dimenB by -\!yM%
  \!yB=3\!dimenA \advance \!yB by -\!dimenB%
  \!yC=2\!dimenB \advance \!yC by -2\!dimenA%
  \!xprime=\!xB  \!yprime=\!yB
  \!dxprime=.5\!xC  \!dyprime=.5\!yC
  \!getf \!midarclength=\!dimenA
  \!getf \advance \!midarclength by 4\!dimenA
  \!getf \advance \!midarclength by \!dimenA
  \divide \!midarclength by 12
  \!arclength=\!dimenA
  \!getf \advance \!arclength by 4\!dimenA
  \!getf \advance \!arclength by \!dimenA
  \divide \!arclength by 12
  \advance \!arclength by \!midarclength
  \global\advance \totalarclength by \!arclength
  \ifdim\!distacross>\!arclength 
    \advance \!distacross by -\!arclength
  \else
    \!initinverseinterp
    \loop\ifdim\!distacross<\!arclength
      \!inverseinterp
      \!xpos=\!t\!xC \advance\!xpos by \!xB
        \!xpos=\!t\!xpos \advance \!xpos by \!xS
      \!ypos=\!t\!yC \advance\!ypos by \!yB
        \!ypos=\!t\!ypos \advance \!ypos by \!yS
      \!plotifinbounds
      \advance\!distacross \plotsymbolspacing
      \!advancedashing
    \repeat  
    \advance \!distacross by -\!arclength
  \fi
  \!xS=\!xE
  \!yS=\!yE
  \ignorespaces}

\def\!getf{\!Pythag\!xprime\!yprime\!dimenA%
  \advance\!xprime by \!dxprime
  \advance\!yprime by \!dyprime}

\def\!initinverseinterp{%
  \ifdim\!arclength>\!zpt
    \!divide{8\!midarclength}\!arclength\!dimenE
    \ifdim\!dimenE<\!wmin \!setinverselinear
    \else 
      \ifdim\!dimenE>\!wmax \!setinverselinear
      \else
        \def\!inverseinterp{\!inversequad}\ignorespaces
         \!removept\!dimenE\!Ew
         \!dimenF=-\!Ew\!dimenE
         \advance\!dimenF by 32pt
         \!dimenG=8pt 
         \advance\!dimenG by -\!dimenE
         \!dimenG=\!Ew\!dimenG
         \!divide\!dimenF\!dimenG\!beta
         \!gamma=1pt
         \advance \!gamma by -\!beta
      \fi
    \fi
  \fi
  \ignorespaces}

\def\!inversequad{%
  \!divide\!distacross\!arclength\!dimenG
  \!removept\!dimenG\!v
  \!dimenG=\!v\!gamma
  \advance\!dimenG by \!beta
  \!dimenG=\!v\!dimenG
  \!removept\!dimenG\!t}

\def\!setinverselinear{%
  \def\!inverseinterp{\!inverselinear}%
  \divide\!dimenE by 8 \!removept\!dimenE\!t
  \!countC=\!intervalno \multiply \!countC 2
  \!countB=\!countC     \advance \!countB -1
  \!countA=\!countB     \advance \!countA -1
  \wlog{\the\!countB th point (\!xmidpt,\!ymidpt) being plotted 
    doesn't lie in the}%
  \wlog{ middle third of the arc between the \the\!countA th 
    and \the\!countC th points:}%
  \wlog{ [arc length \the\!countA\space to \the\!countB]/[arc length 
    \the \!countA\space to \the\!countC]=\!t.}%
  \ignorespaces}
 
\def\!inverselinear{%
  \!divide\!distacross\!arclength\!dimenG
  \!removept\!dimenG\!t}

\def\startrotation{%
  \let\!rotateaboutpivot=\!!rotateaboutpivot
  \let\!rotateonly=\!!rotateonly
  \!ifnextchar{b}{\!getsincos }%
    {\!getsincos by {\!cosrotationangle} {\!sinrotationangle} }}
\def\!getsincos by #1 #2 {%
  \edef\!cosrotationangle{#1}%
  \edef\!sinrotationangle{#2}%
  \!ifcoordmode 
    \let\!ROnext=\!ccheckforpivot
  \else
    \let\!ROnext=\!dcheckforpivot
  \fi
  \!ROnext}
\def\!ccheckforpivot{%
  \!ifnextchar{a}{\!cgetpivot}%
    {\!cgetpivot about {\!xpivotcoord} {\!ypivotcoord} }}
\def\!cgetpivot about #1 #2 {%
  \edef\!xpivotcoord{#1}%
  \edef\!ypivotcoord{#2}%
  \!xpivot=#1\!xunit  \!ypivot=#2\!yunit
  \ignorespaces}
\def\!dcheckforpivot{%
  \!ifnextchar{a}{\!dgetpivot}{\ignorespaces}}
\def\!dgetpivot about #1 #2 {%
  \!xpivot=#1\relax  \!ypivot=#2\relax
  \ignorespaces}

\def\stoprotation{%
  \let\!rotateaboutpivot=\!!!rotateaboutpivot
  \let\!rotateonly=\!!!rotateonly
  \ignorespaces}
 
\def\!!rotateaboutpivot#1#2{%
  \!dimenA=#1\relax  \advance\!dimenA -\!xpivot
  \!dimenB=#2\relax  \advance\!dimenB -\!ypivot
  \!dimenC=\!cosrotationangle\!dimenA
    \advance \!dimenC -\!sinrotationangle\!dimenB
  \!dimenD=\!cosrotationangle\!dimenB
    \advance \!dimenD  \!sinrotationangle\!dimenA
  \advance\!dimenC \!xpivot  \advance\!dimenD \!ypivot
  #1=\!dimenC  #2=\!dimenD
  \ignorespaces}

\def\!!rotateonly#1#2{%
  \!dimenA=#1\relax  \!dimenB=#2\relax 
  \!dimenC=\!cosrotationangle\!dimenA
    \advance \!dimenC -\!rotsign\!sinrotationangle\!dimenB
  \!dimenD=\!cosrotationangle\!dimenB
    \advance \!dimenD  \!rotsign\!sinrotationangle\!dimenA
  #1=\!dimenC  #2=\!dimenD
  \ignorespaces}
\def\!rotsign{}
\def\!!!rotateaboutpivot#1#2{\relax}
\def\!!!rotateonly#1#2{\relax}
\stoprotation

\def\!reverserotateonly#1#2{%
  \def\!rotsign{-}%
  \!rotateonly{#1}{#2}%
  \def\!rotsign{}%
  \ignorespaces}

\def\!getspan span <#1>{%
  \!dshade=#1\relax
  \!ifcoordmode 
    \let\!GRnext=\!GRccheckforAP
  \else
    \let\!GRnext=\!GRdcheckforAP
  \fi
  \!GRnext}
\def\!GRccheckforAP{%
  \!ifnextchar{p}{\!cgetanchor }
    {\!cgetanchor point at {\!xshadesave} {\!yshadesave} }}
\def\!cgetanchor point at #1 #2 {%
  \edef\!xshadesave{#1}\edef\!yshadesave{#2}%
  \!xshade=\!xshadesave\!xunit  \!yshade=\!yshadesave\!yunit
  \ignorespaces}
\def\!GRdcheckforAP{%
  \!ifnextchar{p}{\!dgetanchor}%
    {\ignorespaces}}
\def\!dgetanchor point at #1 #2 {%
  \!xshade=#1\relax  \!yshade=#2\relax
  \ignorespaces}

\def\setshadesymbol{%
  \!ifnextchar<{\!setshadesymbol}{\!setshadesymbol<,,,> }}

\def\!setshadesymbol <#1,#2,#3,#4> (#5#6){%
  \!setputobject{#5}{#6}%
  \setbox\!shadesymbol=\box\!putobject%
  \!shadesymbolxshift=\!xshift \!shadesymbolyshift=\!yshift
  \!dimenA=\!xshift \advance\!dimenA \!smidge
  \!override\!dimenA{#1}\!lshrinkage%
  \!dimenA=\!wd \advance \!dimenA -\!xshift
    \advance\!dimenA \!smidge
    \!override\!dimenA{#2}\!rshrinkage
  \!dimenA=\!dp \advance \!dimenA \!yshift
    \advance\!dimenA \!smidge
    \!override\!dimenA{#3}\!bshrinkage
  \!dimenA=\!ht \advance \!dimenA -\!yshift
    \advance\!dimenA \!smidge
    \!override\!dimenA{#4}\!tshrinkage
  \ignorespaces}
\def\!smidge{-.2pt}%

\def\!override#1#2#3{%
  \edef\!!override{#2}%
  \ifx \!!override\empty
    #3=#1\relax
  \else
    \if z\!!override
      #3=\!zpt
    \else
      \ifx \!!override\!blankz
        #3=\!zpt
      \else
        #3=#2\relax
      \fi
    \fi
  \fi
  \ignorespaces}
\def\!blankz{ z}

\setshadesymbol ({\fiverm .})

\def\!startvshade#1(#2,#3,#4){%
  \let\!!xunit=\!xunit%
  \let\!!yunit=\!yunit%
  \let\!!xshade=\!xshade%
  \let\!!yshade=\!yshade%
  \def\!getshrinkages{\!vgetshrinkages}%
  \let\!setshadelocation=\!vsetshadelocation%
  \!xS=\!M{#2}\!!xunit
  \!ybS=\!M{#3}\!!yunit
  \!ytS=\!M{#4}\!!yunit
  \!shadexorigin=\!xorigin  \advance \!shadexorigin \!shadesymbolxshift
  \!shadeyorigin=\!yorigin  \advance \!shadeyorigin \!shadesymbolyshift
  \ignorespaces}
 
\def\!starthshade#1(#2,#3,#4){%
  \let\!!xunit=\!yunit%
  \let\!!yunit=\!xunit%
  \let\!!xshade=\!yshade%
  \let\!!yshade=\!xshade%
  \def\!getshrinkages{\!hgetshrinkages}%
  \let\!setshadelocation=\!hsetshadelocation%
  \!xS=\!M{#2}\!!xunit
  \!ybS=\!M{#3}\!!yunit
  \!ytS=\!M{#4}\!!yunit
  \!shadexorigin=\!xorigin  \advance \!shadexorigin \!shadesymbolxshift
  \!shadeyorigin=\!yorigin  \advance \!shadeyorigin \!shadesymbolyshift
  \ignorespaces}

\def\!lattice#1#2#3#4#5{%
  \!dimenA=#1
  \!dimenB=#2
  \!countB=\!dimenB
  \!dimenC=#3
  \advance\!dimenC -\!dimenA
  \!countA=\!dimenC
  \divide\!countA \!countB
  \ifdim\!dimenC>\!zpt
    \!dimenD=\!countA\!dimenB
    \ifdim\!dimenD<\!dimenC
      \advance\!countA 1 
    \fi
  \fi
  \!dimenC=\!countA\!dimenB
    \advance\!dimenC \!dimenA
  #4=\!countA
  #5=\!dimenC
  \ignorespaces}

\def\!qshade#1(#2,#3,#4)#5(#6,#7,#8){%
  \!xM=\!M{#2}\!!xunit
  \!ybM=\!M{#3}\!!yunit
  \!ytM=\!M{#4}\!!yunit
  \!xE=\!M{#6}\!!xunit
  \!ybE=\!M{#7}\!!yunit
  \!ytE=\!M{#8}\!!yunit
  \!getcoeffs\!xS\!ybS\!xM\!ybM\!xE\!ybE\!ybB\!ybC
  \!getcoeffs\!xS\!ytS\!xM\!ytM\!xE\!ytE\!ytB\!ytC
  \def\!getylimits{\!qgetylimits}%
  \!shade{#1}\ignorespaces}
 
\def\!lshade#1(#2,#3,#4){%
  \!xE=\!M{#2}\!!xunit
  \!ybE=\!M{#3}\!!yunit
  \!ytE=\!M{#4}\!!yunit
  \!dimenE=\!xE  \advance \!dimenE -\!xS
  \!dimenC=\!ytE \advance \!dimenC -\!ytS
  \!divide\!dimenC\!dimenE\!ytB
  \!dimenC=\!ybE \advance \!dimenC -\!ybS
  \!divide\!dimenC\!dimenE\!ybB
  \def\!getylimits{\!lgetylimits}%
  \!shade{#1}\ignorespaces}
 
\def\!getcoeffs#1#2#3#4#5#6#7#8{%
  \!dimenC=#4\advance \!dimenC -#2
  \!dimenE=#3\advance \!dimenE -#1
  \!divide\!dimenC\!dimenE\!dimenF
  \!dimenC=#6\advance \!dimenC -#4
  \!dimenH=#5\advance \!dimenH -#3
  \!divide\!dimenC\!dimenH\!dimenG
  \advance\!dimenG -\!dimenF
  \advance \!dimenH \!dimenE
  \!divide\!dimenG\!dimenH#8
  \!removept#8\!t
  #7=-\!t\!dimenE
  \advance #7\!dimenF
  \ignorespaces}

\def\!shade#1{%
  \!getshrinkages#1<,,,>\!nil
  \advance \!dimenE \!xS
  \!lattice\!!xshade\!dshade\!dimenE
    \!parity\!xpos
  \!dimenF=-\!dimenF
    \advance\!dimenF \!xE
  \!loop\!not{\ifdim\!xpos>\!dimenF}
    \!shadecolumn%
    \advance\!xpos \!dshade
    \advance\!parity 1
  \repeat
  \!xS=\!xE
  \!ybS=\!ybE
  \!ytS=\!ytE
  \ignorespaces}

\def\!vgetshrinkages#1<#2,#3,#4,#5>#6\!nil{%
  \!override\!lshrinkage{#2}\!dimenE
  \!override\!rshrinkage{#3}\!dimenF
  \!override\!bshrinkage{#4}\!dimenG
  \!override\!tshrinkage{#5}\!dimenH
  \ignorespaces}
\def\!hgetshrinkages#1<#2,#3,#4,#5>#6\!nil{%
  \!override\!lshrinkage{#2}\!dimenG
  \!override\!rshrinkage{#3}\!dimenH
  \!override\!bshrinkage{#4}\!dimenE
  \!override\!tshrinkage{#5}\!dimenF
  \ignorespaces}

\def\!shadecolumn{%
  \!dxpos=\!xpos
  \advance\!dxpos -\!xS
  \!removept\!dxpos\!dx
  \!getylimits
  \advance\!ytpos -\!dimenH
  \advance\!ybpos \!dimenG
  \!yloc=\!!yshade
  \ifodd\!parity 
     \advance\!yloc \!dshade
  \fi
  \!lattice\!yloc{2\!dshade}\!ybpos%
    \!countA\!ypos
  \!dimenA=-\!shadexorigin \advance \!dimenA \!xpos
  \loop\!not{\ifdim\!ypos>\!ytpos}
    \!setshadelocation
    \!rotateaboutpivot\!xloc\!yloc%
    \!dimenA=-\!shadexorigin \advance \!dimenA \!xloc
    \!dimenB=-\!shadeyorigin \advance \!dimenB \!yloc
    \kern\!dimenA \raise\!dimenB\copy\!shadesymbol \kern-\!dimenA
    \advance\!ypos 2\!dshade
  \repeat
  \ignorespaces}
 
\def\!qgetylimits{%
  \!dimenA=\!dx\!ytC              
  \advance\!dimenA \!ytB
  \!ytpos=\!dx\!dimenA
  \advance\!ytpos \!ytS
  \!dimenA=\!dx\!ybC              
  \advance\!dimenA \!ybB
  \!ybpos=\!dx\!dimenA
  \advance\!ybpos \!ybS}
 
\def\!lgetylimits{%
  \!ytpos=\!dx\!ytB
  \advance\!ytpos \!ytS
  \!ybpos=\!dx\!ybB
  \advance\!ybpos \!ybS}
 
\def\!vsetshadelocation{
  \!xloc=\!xpos
  \!yloc=\!ypos}
\def\!hsetshadelocation{
  \!xloc=\!ypos
  \!yloc=\!xpos}

\def\!axisticks {%
  \def\!nextkeyword##1 {%
    \expandafter\ifx\csname !ticks##1\endcsname \relax
      \def\!next{\!fixkeyword{##1}}%
    \else
      \def\!next{\csname !ticks##1\endcsname}%
    \fi
    \!next}%
  \!axissetup
    \def\!axissetup{\relax}%
  \edef\!ticksinoutsign{\!ticksinoutSign}%
  \!ticklength=\longticklength
  \!tickwidth=\linethickness
  \!gridlinestatus
  \!setticktransform
  \!maketick
  \!tickcase=0
  \def\!LTlist{}%
  \!nextkeyword}

\def\ticksout{%
  \def\!ticksinoutSign{+}}

\ticksout

\def\nogridlines{%
  \def\!gridlinestatus{\!gridlinestoofalse}}
\nogridlines

\def\loggedticks{%
  \def\!setticktransform{\let\!ticktransform=\!logten}}
\def\unloggedticks{%
  \def\!setticktransform{\let\!ticktransform=\!donothing}}
\def\!donothing#1#2{\def#2{#1}}
\unloggedticks

\expandafter\def\csname !ticks/\endcsname{%
  \!not {\ifx \!LTlist\empty}
    \!placetickvalues
  \fi
  \def\!tickvalueslist{}%
  \def\!LTlist{}%
  \expandafter\csname !axis/\endcsname}

\def\!maketick{%
  \setbox\!boxA=\hbox{%
    \beginpicture
      \!setdimenmode
      \setcoordinatesystem point at {\!zpt} {\!zpt}   
      \linethickness=\!tickwidth
      \ifdim\!ticklength>\!zpt
        \putrule from {\!zpt} {\!zpt} to
          {\!ticksinoutsign\!tickxsign\!ticklength}
          {\!ticksinoutsign\!tickysign\!ticklength}
      \fi
      \if!gridlinestoo
        \putrule from {\!zpt} {\!zpt} to
          {-\!tickxsign\!xaxislength} {-\!tickysign\!yaxislength}
      \fi
    \endpicturesave <\!Xsave,\!Ysave>}%
    \wd\!boxA=\!zpt}
  
\def\!ticksin{%
  \def\!ticksinoutsign{-}%
  \!maketick
  \!nextkeyword}

\def\!ticksout{%
  \def\!ticksinoutsign{+}%
  \!maketick
  \!nextkeyword}

\def\!tickslength<#1> {%
  \!ticklength=#1\relax
  \!maketick
  \!nextkeyword}

\def\!tickslong{%
  \!tickslength<\longticklength> }

\def\!ticksshort{%
  \!tickslength<\shortticklength> }

\def\!tickswidth<#1> {%
  \!tickwidth=#1\relax
  \!maketick
  \!nextkeyword}

\def\!ticksandacross{%
  \!gridlinestootrue
  \!maketick
  \!nextkeyword}

\def\!ticksbutnotacross{%
  \!gridlinestoofalse
  \!maketick
  \!nextkeyword}

\def\!tickslogged{%
  \let\!ticktransform=\!logten
  \!nextkeyword}

\def\!ticksunlogged{%
  \let\!ticktransform=\!donothing
  \!nextkeyword}

\def\!ticksunlabeled{%
  \!tickcase=0
  \!nextkeyword}

\def\!ticksnumbered{%
  \!tickcase=1
  \!nextkeyword}

\def\!tickswithvalues#1/ {%
  \edef\!tickvalueslist{#1! /}%
  \!tickcase=2
  \!nextkeyword}

\def\!ticksquantity#1 {%
  \ifnum #1>1
    \!updatetickoffset
    \!countA=#1\relax
    \advance \!countA -1
    \!ticklocationincr=\!axisLength
      \divide \!ticklocationincr \!countA
    \!ticklocation=\!axisstart
    \loop \!not{\ifdim \!ticklocation>\!axisend}
      \!placetick\!ticklocation
      \ifcase\!tickcase
          \relax 
        \or
          \relax 
        \or
          \expandafter\!gettickvaluefrom\!tickvalueslist
          \edef\!tickfield{{\the\!ticklocation}{\!value}}%
          \expandafter\!listaddon\expandafter{\!tickfield}\!LTlist%
      \fi
      \advance \!ticklocation \!ticklocationincr
    \repeat
  \fi
  \!nextkeyword}

\def\!ticksat#1 {%
  \!updatetickoffset
  \edef\!Loc{#1}%
  \if /\!Loc
    \def\next{\!nextkeyword}%
  \else
    \!ticksincommon
    \def\next{\!ticksat}%
  \fi
  \next}    
      
\def\!ticksfrom#1 to #2 by #3 {%
  \!updatetickoffset
  \edef\!arg{#3}%
  \expandafter\!separate\!arg\!nil
  \!scalefactor=1
  \expandafter\!countfigures\!arg/
  \edef\!arg{#1}%
  \!scaleup\!arg by\!scalefactor to\!countE
  \edef\!arg{#2}%
  \!scaleup\!arg by\!scalefactor to\!countF
  \edef\!arg{#3}%
  \!scaleup\!arg by\!scalefactor to\!countG
  \loop \!not{\ifnum\!countE>\!countF}
    \ifnum\!scalefactor=1
      \edef\!Loc{\the\!countE}%
    \else
      \!scaledown\!countE by\!scalefactor to\!Loc
    \fi
    \!ticksincommon
    \advance \!countE \!countG
  \repeat
  \!nextkeyword}

\def\!updatetickoffset{%
  \!dimenA=\!ticksinoutsign\!ticklength
  \ifdim \!dimenA>\!offset
    \!offset=\!dimenA
  \fi}

\def\!placetick#1{%
  \if!xswitch
    \!xpos=#1\relax
    \!ypos=\!axisylevel
  \else
    \!xpos=\!axisxlevel
    \!ypos=#1\relax
  \fi
  \advance\!xpos \!Xsave
  \advance\!ypos \!Ysave
  \kern\!xpos\raise\!ypos\copy\!boxA\kern-\!xpos
  \ignorespaces}

\def\!gettickvaluefrom#1 #2 /{%
  \edef\!value{#1}%
  \edef\!tickvalueslist{#2 /}%
  \ifx \!tickvalueslist\!endtickvaluelist
    \!tickcase=0
  \fi}
\def\!endtickvaluelist{! /}

\def\!ticksincommon{%
  \!ticktransform\!Loc\!t
  \!ticklocation=\!t\!!unit
  \advance\!ticklocation -\!!origin
  \!placetick\!ticklocation
  \ifcase\!tickcase
    \relax 
  \or 
    \ifdim\!ticklocation<-\!!origin
      \edef\!Loc{$\!Loc$}%
    \fi
    \edef\!tickfield{{\the\!ticklocation}{\!Loc}}%
    \expandafter\!listaddon\expandafter{\!tickfield}\!LTlist%
  \or 
    \expandafter\!gettickvaluefrom\!tickvalueslist
    \edef\!tickfield{{\the\!ticklocation}{\!value}}%
    \expandafter\!listaddon\expandafter{\!tickfield}\!LTlist%
  \fi}

\def\!separate#1\!nil{%
  \!ifnextchar{-}{\!!separate}{\!!!separate}#1\!nil}
\def\!!separate-#1\!nil{%
  \def\!sign{-}%
  \!!!!separate#1..\!nil}
\def\!!!separate#1\!nil{%
  \def\!sign{+}%
  \!!!!separate#1..\!nil}
\def\!!!!separate#1.#2.#3\!nil{%
  \def\!arg{#1}%
  \ifx\!arg\!empty
    \!countA=0
  \else
    \!countA=\!arg
  \fi
  \def\!arg{#2}%
  \ifx\!arg\!empty
    \!countB=0
  \else
    \!countB=\!arg
  \fi}
 
\def\!countfigures#1{%
  \if #1/%
    \def\!next{\ignorespaces}%
  \else
    \multiply\!scalefactor 10
    \def\!next{\!countfigures}%
  \fi
  \!next}

\def\!scaleup#1by#2to#3{%
  \expandafter\!separate#1\!nil
  \multiply\!countA #2\relax
  \advance\!countA \!countB
  \if -\!sign
    \!countA=-\!countA
  \fi
  #3=\!countA
  \ignorespaces}

\def\!scaledown#1by#2to#3{%
  \!countA=#1\relax
  \ifnum \!countA<0 
    \def\!sign{-}
    \!countA=-\!countA
  \else
    \def\!sign{}%
  \fi
  \!countB=\!countA
  \divide\!countB #2\relax
  \!countC=\!countB
    \multiply\!countC #2\relax
  \advance \!countA -\!countC
  \edef#3{\!sign\the\!countB.}
  \!countC=\!countA 
  \ifnum\!countC=0 
    \!countC=1
  \fi
  \multiply\!countC 10
  \!loop \ifnum #2>\!countC
    \edef#3{#3\!zero}%
    \multiply\!countC 10
  \repeat
  \edef#3{#3\the\!countA}
  \ignorespaces}

\def\!placetickvalues{%
  \advance\!offset \tickstovaluesleading
  \if!xswitch
    \setbox\!boxA=\hbox{%
      \def\\##1##2{%
        \!dimenput {##2} [B] (##1,\!axisylevel)}%
      \beginpicture 
        \!LTlist
      \endpicturesave <\!Xsave,\!Ysave>}%
    \!dimenA=\!axisylevel
      \advance\!dimenA -\!Ysave
      \advance\!dimenA \!tickysign\!offset
      \if -\!tickysign
        \advance\!dimenA -\ht\!boxA
      \else
        \advance\!dimenA  \dp\!boxA
      \fi
    \advance\!offset \ht\!boxA 
      \advance\!offset \dp\!boxA
    \!dimenput {\box\!boxA} [Bl] <\!Xsave,\!Ysave> (\!zpt,\!dimenA)
  \else
    \setbox\!boxA=\hbox{%
      \def\\##1##2{%
        \!dimenput {##2} [r] (\!axisxlevel,##1)}%
      \beginpicture 
        \!LTlist
      \endpicturesave <\!Xsave,\!Ysave>}%
    \!dimenA=\!axisxlevel
      \advance\!dimenA -\!Xsave
      \advance\!dimenA \!tickxsign\!offset
      \if -\!tickxsign
        \advance\!dimenA -\wd\!boxA
      \fi
    \advance\!offset \wd\!boxA
    \!dimenput {\box\!boxA} [Bl] <\!Xsave,\!Ysave> (\!dimenA,\!zpt)
  \fi}

\normalgraphs
\catcode`!=12 

\catcode`@=11 \catcode`!=11
  
\let\!pictexendpicture=\endpicture 
\let\!pictexframe=\frame
\let\!pictexlinethickness=\linethickness
\let\!pictexmultiput=\multiput
\let\!pictexput=\put

\def\beginpicture{%
  \setbox\!picbox=\hbox\bgroup%
  \let\endpicture=\!pictexendpicture
  \let\frame=\!pictexframe
  \let\linethickness=\!pictexlinethickness
  \let\multiput=\!pictexmultiput
  \let\put=\!pictexput
  \let\input=\@@input   
  \!xleft=\maxdimen  
  \!xright=-\maxdimen
  \!ybot=\maxdimen
  \!ytop=-\maxdimen}

\let\frame=\!latexframe

\let\pictexframe=\!pictexframe

\let\linethickness=\!latexlinethickness
\let\pictexlinethickness=\!pictexlinethickness

\let\\=\@normalcr
\catcode`@=12 \catcode`!=12

\catcode`\@=11

\@input{picmore.tex}

\@input{piccorr.sty}

\setlength{\hoffset}{-1.cm}
\setlength{\voffset}{-1.cm}
\setlength{\textwidth}{ 14.8cm}
\setlength{\textheight}{22cm}
\setlength{\parindent}{8mm}
\frenchspacing
\font\teneufm=eufm10
\font\seveneufm=eufm7
\font\fiveeufm=eufm5
\newfam\eufmfam
\textfont\eufmfam=\teneufm
\scriptfont\eufmfam=\seveneufm
\scriptscriptfont\eufmfam=\fiveeufm

\def\hh{{\rm h}}

\def\mm{{\rm m}}
\def\MM{{\rm M}}
\def\rr{{\rm r}}

\def\UU{{\rm U}}
\newcommand\beq[1]{ \begin{equation}\label{#1} }
\newcommand{\eeq}{ \end{equation} }
\newcommand{\beqno}{ \[ }
\newcommand{\eeqno}{ \] }
\newcommand\beqa[1]{ \begin{eqnarray} \label{#1}}
\newcommand{\eeqa}{ \end{eqnarray} }
\newcommand{\beqano}{ \begin{eqnarray*} }
\newcommand{\eeqano}{ \end{eqnarray*} }
\newcommand\arr[1]{\left\{\begin{array}{l}#1\end{array}\right.}
\renewcommand{\theequation}{\arabic{section}.\arabic{equation}}

\newtheorem{theorem}{Theorem}[section]
\newtheorem{definition}{Definition}[section]
\newtheorem{proposition}{Proposition}[section]
\newtheorem{lemma}{Lemma}[section]
\newtheorem{sublemma}{Sublemma}[section]
\newtheorem{remark}{Remark}[section]
\newtheorem{notationalremark}{Notations}[section]
\newtheorem{corollary}{Corollary}[section]
\newtheorem{assumption}{Assumption}[section]
\newtheorem{claim}{Claim}[section]
\newtheorem{tools}{$\negsp\negsp$}[subsection]

\newcommand\thm[1]{ \begin{theorem}\label{#1}}
\newcommand\thmtwo[2]{ \begin{theorem}[#1]\label{#2}}
\newcommand\ethm{ \end{theorem} }
\newcommand\dfn[1]{ \begin{definition}\label{#1} \rm}
\newcommand\dfntwo[2]{ \begin{definition}[#1]\label{#2} \rm}
\newcommand\edfn{ \end{definition} }
\newcommand\pro[1]{ \begin{proposition}\label{#1}}
\newcommand\protwo[2]{ \begin{proposition}[#1]\label{#2}}
\newcommand\epro{ \end{proposition} }
\newcommand\lem[1]{ \begin{lemma}\label{#1}}
\newcommand\lemtwo[2]{ \begin{lemma}[#1]\label{#2}}
\newcommand\elem{ \end{lemma} }
\newcommand\sublem[1]{ \begin{sublemma}\label{#1}}
\newcommand\sublemtwo[2]{ \begin{sublemma}[#1]\label{#2}}
\newcommand\esublem{ \end{sublemma} }
\newcommand\rem[1]{ \begin{remark}\label{#1} \rm}
\newcommand\erem{ \end{remark} }
\newcommand\notrem[1]{ \begin{notationalremark}\label{#1} \rm}
\newcommand\enotrem{ \end{notationalremark} }
\newcommand\cor[1]{ \begin{corollary}\label{#1}}
\newcommand\cortwo[2]{ \begin{corollary}[#1]\label{#2}}
\newcommand\ecor{ \end{corollary} }
\newcommand\asmp[1]{ \begin{assumption}\label{#1}}
\newcommand\asmptwo[2]{ \begin{assumption}[#1]\label{#2}}
\newcommand\easmp{ \end{assumption} }
\newcommand\clm[1]{ \begin{claim}\label{#1}}
\newcommand\eclm{ \end{claim} }
\newcommand{\proof}{\par\medskip\noindent{\bf Proof\ }}
\newcommand\equ[1]{{\rm (\ref{#1})}}

\expandafter\chardef\csname pre amssym.def
at\endcsname=\the\catcode`\@
\catcode`\@=11
\def\undefine#1{\let#1\undefined}
\def\newsymbol#1#2#3#4#5{\let\next@\relax
 \ifnum#2=\@ne\let\next@\msafam@\else
 \ifnum#2=\tw@\let\next@\msbfam@\fi\fi
 \mathchardef#1="#3\next@#4#5}
\def\mathhexbox@#1#2#3{\relax
 \ifmmode\mathpalette{}{\m@th\mathchar"#1#2#3}%
 \else\leavevmode\hbox{$\m@th\mathchar"#1#2#3$}\fi}
\def\hexnumber@#1{\ifcase#1 0\or 1\or 2\or 3\or 4\or 5\or 6\or 7\or
8\or
 9\or A\or B\or C\or D\or E\or F\fi}
\ifcase\@ptsize
 \font\tenmsb=msbm10
 \font\sevenmsb=msbm7
 \font\fivemsb=msbm5
\or
 \font\tenmsb=msbm10 scaled \magstephalf
 \font\sevenmsb=msbm7 scaled \magstephalf
 \font\fivemsb=msbm5  scaled \magstephalf
\or
 \font\tenmsb=msbm10 scaled \magstep1
 \font\sevenmsb=msbm7 scaled \magstep1
 \font\fivemsb=msbm5 scaled \magstep1
\fi
\newfam\msbfam
\textfont\msbfam=\tenmsb
\scriptfont\msbfam=\sevenmsb
\scriptscriptfont\msbfam=\fivemsb
\edef\msbfam@{\hexnumber@\msbfam}
\def\Bbb#1{\fam\msbfam\relax#1}
\def\widehat#1{\setboxz@h{$\m@th#1$}%
 \ifdim\wdz@>\tw@ em\mathaccent"0\msbfam@5B{#1}%
 \else\mathaccent"0362{#1}\fi}
\def\widetilde#1{\setboxz@h{$\m@th#1$}%
 \ifdim\wdz@>\tw@ em\mathaccent"0\msbfam@5D{#1}%
 \else\mathaccent"0365{#1}\fi}

\def\RIfM@{\relax\ifmmode}
\def\nonmatherr@#1{\errmessage{\string#1\space allowed only in math mode}}
\def\Bbb{\RIfM@\expandafter\Bbb@\else
 \expandafter\nonmatherr@\expandafter\Bbb\fi}
\def\Bbb@#1{{\Bbb@@{#1}}}
\def\Bbb@@#1{\fam\msbfam\relax#1}
\def\setboxz@h{\setbox\z@\hbox}
\def\wdz@{\wd\z@}
\catcode`\@=\csname pre amssym.def at\endcsname
\newcommand{\nl}{{\smallskip\noindent}}

\newcommand{\negsp}{\hspace{-.09truecm}}  

\newcommand{\dst}{\displaystyle}
\newcommand\ovl[1]{ \overline {#1} }

\newcommand{\torus}{ {\Bbb T}   }
\renewcommand{\natural}{ {\Bbb N}   }
\newcommand{\real}{ {\Bbb R}   }
\newcommand{\integer}{ {\Bbb Z}   }
\newcommand{\complex}{ {\Bbb C}   }

\renewcommand{\a }{ {\alpha}   }
\renewcommand{\b}{ {\beta}   }
\newcommand{\g}{ {\gamma}   }

\renewcommand{\d}{ {\delta}   }
\newcommand{\D}{ {\Delta}   }

\newcommand{\HH }{ {\rm H}   }

\renewcommand{\l}{ {\lambda}   }
\renewcommand{\L}{ {\Lambda}   }
\newcommand{\m}{ {\mu}   }
\newcommand{\n}{ {\nu}   }

\newcommand{\p}{ {\pi}   }
\renewcommand{\P}{ {\Pi}   }
\renewcommand{\r}{ {\rho}   }
\newcommand{\s}{ {\sigma}   }

\newcommand{\f}{ {\varphi}   }

\renewcommand{\o}{ {\omega}   }

\newcommand{\GG}{{\rm G}}

\renewcommand{\=}{ {\, :=\, }   }
\newcommand{\cB}{ {\cal B} }
\newcommand{\cE}{ {\cal E} }

\newcommand{\cH}{ {\cal H} }
\newcommand{\cK}{ {\cal K} }

\newcommand{\cD}{ {\cal D} }

\newcommand{{\cG}}{{\cal G} }
\newcommand{\cM}{ {\cal M} }

\newcommand{\cP}{ {\cal P} }
\newcommand{\cI}{ {\cal I} }
\newcommand{\cJ}{ {\cal J} }

\newcommand{\cS}{ {\cal S} }

\newcommand\RR{{\rm R}}
\newcommand\CC{{\rm C}}
\newcommand\ee{{\rm e}}
\newcommand\ZZ{{\rm Z}}
\newcommand\zz{{\rm z}}

\begin{document}

\title{An analysis of the  Sun--Earth--Asteroid systems based on the two--centre problem\footnote{ 
gravitation.
}}

\author{  
Gabriella Pinzari\thanks{The author is partially supported by the H2020 Project ERC Starting Grant 677793 StableChaoticPlanetM. }\\
\footnotesize{Dipartimento di Matematica ``Tullio Levi Civita''} \\
\vspace{-.2truecm}
\footnotesize{Universit\`a di Padova}\\
{\scriptsize gabriella.pinzari@math.unipd.it}
}\date{May 26, 2017}
\maketitle

\begin{abstract}\footnotesize{
We propose a new analysis of the two--centre problem particularly suited to be used as a basis to study  the dynamics of  Sun--Earth--Asteroid systems. Our method, based on a tricky choice of initial coordinates, allows us to  evaluate the risk of collisions between the Asteroid and Earth. Moreover, it seems to be fitted to prove the existence of chains of transition tori in the planar  Sun--Earth--Asteroid systems.
}
\end{abstract}

\maketitle

\newpage
\tableofcontents
\newpage\section{Introduction}
\setcounter{equation}{0}
\renewcommand{\theequation}{\arabic{equation}}
The Law of Universal Gravitation, according to which, any two masses in the Universe attract each other with a law going as the inverse squared distance, was stated in 1687  by   Isaac Newton. Newton  was aimed  to find a  theoretical explanation to the laws discovered by Johannes Kepler between 1609 and 1619. 
At the same time, he provided the  {\it exact solution} of the simplest
gravitational system: the
two--body problem (2{\sc bp}), or: the problem of Sun and Earth. He also tried to attack  the analogous problem with three masses (3{\sc bp}: Sun, Earth and Moon), and  then gave up,  
calling it a  `head ache problem'. In 1899 
 Henri Poincar\'e proved the {\it non--integrability} of 3{\sc bp}, and 
 this motivated him to  
introduce the 
 concept of {\it chaos} in mathematics, \cite{MR0087813}.  A  major breakthrough 
came from   {\it Kolmogorov--Arnold--Moser} ({\sc kam}) {\it theory}, developed starting with the fundamental papers \cite{kolmogorov54, moser1962, arnold63c}. Rather than finding an
 explicit time law for the coordinate functions,
{\sc kam} theory turned to have, as main objective, that of providing an estimate of `stable motions', from the probabilistic (measure--theoretic) point of view. Clare answers about the `metric' stability of the {\it planetary problem} (the problem of one star and any number of  smaller masses) have been given in  \cite{arnold63, robutel95, fejoz04, pinzari-th09, chierchiaPi11b}. \\
Much  less known as further example of {exactly solved} gravitational system is the so--called {\it two--centre problem} (2{\sc cp}),  solved by Euler in XVIII century. It consists of one particle, attracted by two {\it fixed} masses. 
The model might, at first sight, seem unsatisfactory from the physical point of view, but, as one might argue, 
may be regarded as  good initial approximation to 3{\sc bp}.
Its solution is  given in the  form of a  non--linear system involving  elliptic integrals \cite[Eq. (53)]{bekovO78}.  Recently, new analysis of the problem have been worked out, by \cite{waalkensDHR04, biscaniI2016}, to which papers we refer  for an overview on motivations and
complete references. Maybe due to the difficulty of handling, at practical level, Euler's solutions,   2{\sc cp} 
has been   not frequently used  in the study of 3{\sc bp}. A study 
in this direction goes back to  \cite{charlier1907}, who applied 
2{\sc cp} 
 to  the {\it restricted}  3{\sc bp}, a  model of the true system that will be recalled below. \\ 
 In this paper, we  present an analysis of 2{\sc cp} that  allows us to write at least a first--order (in the ratio $\varepsilon$ of the masses of the attracting centers) solution {\it explicitly}. Therefore, the method works well when $\varepsilon$ is small.
 The novelty with respect to previous analysis  is that the major information on the dynamics generated by 2{\sc cp} is deferred to the one of its  `Euler integral' (the first integral, found by Euler,  that determines the integrability of 2{\sc cp}). More precisely, we prove that,  in a suitable set of coordinates, 2{\sc cp}'s Hamiltonian and  its Euler integral have the same trajectories, up to a rescaling of time. 
Since, in such coordinates, the Euler integral is much simpler than the original two--centre Hamiltonian, at least at the first order in the masses' ratio, our trick consists in studying the dynamics of the latter function.
The method looks particularly amenable in  the planar problem, because, in this case, the phase portrait of the Euler integral may be studied  exactly, or at least its leading part
 (this is done Section~\ref{G0levels}; see also Figure~\ref{fig: stable}). In  the case of the spatial problem, the procedure can still be applied, but   the analysis would require to solve  a cubic equation (in the planar case, the analogous equation reduces to order two).  However, we do not insist to  develop the theory of the spatial problem (which certainly is the next step of the research). Rather, aiming highlight the utility of the method from the concrete point of view, we 
  discuss an application  to a particular 3{\sc bp}, namely, the planar
 {\it Sun--Earth--Asteroid} system ({\sc sea}).
Most of times, {\sc sea} is studied, from the theoretical point of view, as a {\it restricted three--body problem}. 
 This is a model 
 where the  two most massive  bodies are constrained on circular, co--planar trajectories  having a common centre. A third small body is attracted  by the two, without the two are attracted by it. 
Notwithstanding the important results that have been obtained from the study of the restricted problem, 
the model is affected by an important limitation: the too low (two)   number of degree of freedom
 causes, by {\sc kam} theory, a {\it  confinement at all times} of action coordinates, in the sense of \cite{nehorosev77}. Such confinement is not expected to hold for the true system, since the real Hamiltonian has three, four degrees of freedom in the planar, spatial case, respectively.\\
We  propose an alternative analysis of {\sc sea},  based on its full Hamiltonian, rather than a model.
We write such Hamiltonian as a small perturbation of 2{\sc cp}. 
Using our new approach to 2{\sc cp}, we prove, in the planar {\sc sea}, the existence of stable motions,  with the pericentre
of the Asteroid performing librations or complete rotations,
 even in the case that the orbit of the Asteroid around the Sun encloses Earth -- a situation, geometrically, at risk of collisions between the two. As a byproduct of the proof, the risk of collision may be evaluated simply looking at the numerical  value of the aforementioned Euler integral. Indeed, we find that collisions occur only if (lowest order approximation) of Euler integral takes a suitable value, depending on the mass of Earth and its distance from the sun.

 \nl
 This paper is organized as follows. In Section~\ref{main} we discuss the procedure  (based in an essential way on a good choice of canonical coordinates) that allows us to find the equivalent Hamiltonian, and we discuss its application to {\sc sea}. In Section~\ref{perspectives} we draw conclusions and foresee perspectives of future work. In particular, we conjecture  the existence of chains of {\it transition tori},  in the sense of \cite{arnold64}, in  the planar {\sc sea}. In order to keep the paper as much readable as possible, we relegate the most technical parts to the appendices.

 \newpage
 \section{SEA system via 2CP}\label{main}
\subsection{The two--centre problem}\label{sec: 2C}
The two--centre problem is the problem of determining the motions of one moving  mass $m$ gravitationally attracted by two fixed masses $\MM$, $\MM'$.

\nl
Let us  fix a orthonormal frame $({\mathbf i},{\mathbf j},{\mathbf k})$ in $\real^3$.
After changing the time $t':=g m \MM t$, and  with  $\mm:=g m^2\MM$, $\varepsilon:=\MM'/\MM$, where $g$ is the gravity constant, we write the Hamiltonian  as
\beq{2C}{\rm h}:=\frac{|{\mathbf y}|^2}{2{{\mm}}}-\frac{1}{|{\mathbf x}|}-\frac{\varepsilon  }{|{\mathbf x}'-{\mathbf x}|}\qquad {\mathbf x}\notin\big\{ {\mathbf 0},\ {\mathbf x}'\big\}\eeq
 where  ${\mathbf y}$, ${\mathbf x}$ are impulse--position coordinates of the attracted body (${\mathbf y}=\mm\dot{\mathbf x}$, in the new time), while $|{\mathbf x}|$ denotes Euclidean norm.  Note that the two attracting centers have been posed at ${\mathbf 0}$, ${\mathbf x}'$, rather than, as more commonly done, at symmetric positions (e.g., $\pm {\mathbf i}$) with respect to the origin. For generality, we refer to the spatial problem, namely, ${\mathbf y}$, ${\mathbf x}'$, ${\mathbf x}\in \real^3$. Later on, we shall reduce to the planar case as a sub--case  of the spatial one.

\subsubsection{First integrals}
In view of the application to {\sc sea}, we regard $\hh$ as a {\it six--degrees of freedom system}, i.e., as a function of $({\mathbf y}',{\mathbf y},{\mathbf x}',{\mathbf x})\in \real^3\times\real^3\times\real^3\times\real^3\setminus\{{\mathbf x}=0\ ,\ {\mathbf x}={\mathbf x}'\} $, which is constant with respect to ${\mathbf y}'$. 
In such a enlarged phase space, 
the vectors  ${\mathbf x}'$ and   ${\mathbf C}_{\rm t}={\mathbf C}'+{\mathbf C}$ (the total angular momentum vector), where ${\mathbf C}'={\mathbf x}'\times {\mathbf y}'$, ${\mathbf C}={\mathbf x}\times {\mathbf y}$, are first integrals to $\hh$. We then have {\it six} conserved quantities, which however are not all mutually Poisson--commuting. It is nevertheless possible to extract, out of such six quantities,  the following four commuting ones:
\beq{trivial integrals}
 {\rm Z}:={\mathbf C}_{\rm t}\cdot {\mathbf k}\ ,\  \dst \GG:=|{\mathbf C}_{\rm t}|\ ,\ \dst \Theta:=\frac{{\mathbf C}\cdot {\mathbf x}'}{|{\mathbf x}'|}\ ,\ \rr'=|{\mathbf x}'|\ .\eeq
Obviously, one expects that  two more  first integrals, independent of~\equ{trivial integrals}, can be chosen. This aspect will be discussed in Section~\ref{coordinates} below.

\nl
The  integrability of $\hh$  relies on the existence of a further independent, commuting first integral, found by Euler. It is defined  as follows. Let $\ee$, ${\mathbf P}$, with $\ee\ne 0$ and $|{\mathbf P}|=1$, be the eccentricity and the pericentre direction of the Keplerian orbit associated to 
 \beq{one centre}\cK=\frac{|{\mathbf y}|^2}{2\mm}-\frac{1}{|{\mathbf x}|}\ .\eeq 
Then Euler's first integral to~\equ{2C} is 
\beq{cal G}
{\cal G}={\cal G}_0+\varepsilon{\cal G}_1
\eeq
where
 \beq{G0}{\cal G}_{0}:=
|{\mathbf C}|^2-\mm{\rm e}\, {\mathbf x}' \cdot{\mathbf P}\ ,\ \cG_{1}:=\mm\frac{({\mathbf x}'-{\mathbf x})\cdot {\mathbf x}'}{|{\mathbf x}'-{\mathbf x}|}\ .\eeq
The formula~\equ{G0}, in this precise form, is not standard in the literature, since usually Euler's integral is written in terms of elliptic coordinates. A derivation of it may be found in Appendix~\ref{App: 2C}.\\
Note that $\cG_0$ is a combination of first integrals to $\cK$, as one should expect, since $\hh$ reduces to $\cK$ when $\varepsilon=0$. This also explains why, as one can see from~\equ{G0}, $\cG_0$ has no singularities for ${\mathbf x}={\mathbf x}'$, while $\cG_1$ and $\hh$ do have.

\subsubsection
{Kepler maps and $\l$--normalization}\label{normalization}
We write  
$\hh=\cK+\varepsilon\UU$
where $\cK$ is the Keplerian term~\equ{one centre}, while $\UU=-\frac{1}{|{\mathbf x}'-{\mathbf x}|}$.
We look at systems of canonical coordinates of the form
$\textrm{\sc k}=\big((\L,u,v),\l\big)\subset  \real^{11}\times \torus$ under which (for $\cK<0$) this term takes the form 
\beq{Kepler}\hh_\cK(\L)=-\frac{\mm}{2\L^2}\eeq 
with $\l$ having the form $\l=\ell+\f(\L,u,v)$, where
\beq{a}\L=\sqrt{\mm a}\ ,\quad \ell=\textrm{\rm mean anomaly of }{\mathbf x}_{\textrm{\sc k}}
\eeq
with $a$ the semi--major axis of the ellipse generated by~\equ{one centre}.
We call such kind of systems {\it partial Kepler maps}.
Then  $\hh$ becomes \beq{hk}\hh_{\textrm{\sc k}}\big((\L,u,v),\l\big)=-\frac{\mm}{2\L^2}-\frac{\varepsilon}{|{\mathbf x}'_{\textrm{\sc k}}(\L,u,v)-{\mathbf x}_{\textrm{\sc k}}\big((\L,u,v),\l\big)|}\ ,\eeq
with ${\mathbf x}_{\textrm{\sc k}}$, ${\mathbf x}_{\textrm{\sc k}}$ denoting ${\mathbf x}'$, ${\mathbf x}$ written in terms of $\textrm{\sc k}$. Note that, using a terminology introduced in  \cite{arnold63},   $\hh_{\textrm{\sc k}}$ is a `properly--degenerate' Hamiltonian,  in the sense that its unperturbed part depends on a number of action coordinates strictly less than the number of degrees of freedom. This is a quite common fact in gravitational problems.

\nl
We have to fix a domain for the coordinates so as to exclude collisions. We impose a strong non--collision condition taking 
$\cD =\D^{\rm c}\times\torus$, where
\beq{B0}\D:=\Big\{(\L, u,v):\ \exists\l:\ {\mathbf x}'_{\textrm{\sc k}}(\L,u,v)= {\mathbf x}_{\textrm{\sc k}}\big((\L,u,v),\l\big)\ \Big\}\ .\eeq
Note, for next need, that $\D$ can be  split as a product 
\beq{D}
\D=\cS\times\P
\eeq where:

\begin{itemize}
\item[--] $\P$ is the set of $(\L,u,v)$ such that ${\mathbf x}'_{\textrm{\sc k}}(\L,u,v)$ lies on the plane of the orbit $\l\to x_{\textrm{\sc k}}$, and hence (being orthogonal to ${\mathbf C}_{\textrm{\sc k}}(\L,u,v)$)
\beq{Pi}
\P=\big\{(\L,u,v):\ \Theta=0\big\}\ ;
\eeq
\item[--]  $\cS$ is defined so that, if $\theta(\L,u,v)$ denotes
 is the convex angle formed by ${\mathbf  x}'_{\textrm{\sc k}}(\L,u,v)$ and ${\mathbf  P}_{\textrm{\sc k}}(\L,u,v)$,
 namely,
  \beq{costheta}\theta:\ \cos\theta=\frac{{\mathbf  x}'_{\textrm{\sc k}}\cdot{\mathbf  P}_{\textrm{\sc k}}}{\rr'}\eeq
  then
 the distance $\rr'(\L,u,v)=|{\mathbf  x}'_{\textrm{\sc k}}|$ of ${\mathbf  x}'_{\textrm{\sc k}}$ from the sun satisfies the equality  \beq{collisions}\rr'=\frac{a(1-\ee^2)}{1+\ee\cos\theta}\ .\eeq
Using~\equ{costheta} into~\equ{collisions} and 
recalling  that $a(1-\ee^2)={\rm G}^2/\mm$ and the definition of $\cG_0$ in~\equ{G0}, 
  one sees that
\beq{S}\cS=\{(\L,u,v):\ {\cal G}_0=\mm\rr'\}\ .\eeq
\end{itemize}

\subsubsection{Streaming ideas}
Before, we switch to details, let us try to give an informal account on the  ideas we put in play in the next sections. 
One reasonably expects that, on the collision--less set $\cD$ defined above, due to the close--to--be integrable and simultaneously integrable properties of $\hh_{\textrm{\sc k}}$ in~\equ{hk},
one can  eliminate,  in $\hh_{\textrm{\sc k}}$, the dependence on $\l$ via a {\it convergent} perturbative procedure (in the parameter $\varepsilon$), so as  to transform, via an $\varepsilon$--close to the identity canonical transformation
\beq{eq: norm}\textrm{\sc k}=(\L,\l, u,v)\to \ovl{\textrm{\sc k}}=(\ovl\L,\ovl\l, \ovl u,\ovl v)\eeq
the Hamiltonian $\hh_{\textrm{\sc k}}$ in~\equ{hk} to
\beq{eq: normalization}\ovl\hh:=\hh_{\ovl{\textrm{\sc k}}}=-\frac{\mm}{2\L^2}+\varepsilon\ovl\UU+\cdots\eeq
where $\ovl\UU$ is the $\ovl\l$--average of $\UU$, etc.
Because of the proper degeneracy mentioned above, there would be, in principle, many ways to obtain~\equ{eq: normalization},
depending on how many choices of {\it partial Kepler maps} coordinates one has at our disposal.
However, whatever is the choice of $\textrm{\sc k}$, one  can prove that {\it the function  \beq{ovl G}\ovl{\cal G}:=\cG_{\ovl{\textrm{\sc k}}}\eeq corresponding to  ${\cal G}$ in~\equ{cal G}, written in the  coordinates at right hand side in~\equ{eq: norm}, is itself $\ovl\l$--independent} (details are in Appendix~\ref{lambda independence}). Now, in the next Section~\ref{coordinates}, we shall present a carefully chosen Kepler map such that $\hh_{\textrm{\sc k}}$ and $\cG_{\textrm{\sc k}}$ have {\it two} effective degrees of freedom (they depend on two angle coordinates only). The trick is that $\textrm{\sc k}$ includes, among its coordinates, many first integrals of $\hh$ (all of them, but $\cG$). In this situation, the corresponding functions $\ovl\hh$, $\ovl\cG$ in~\equ{eq: normalization}--\equ{ovl G} have just {\it one} effective degree of freedom and, in addition, Poisson--commute. This implies that, up to rescaling of time (that will be quantified in Section~\ref{norm}), they have the same trajectories. 
But since $\cG_0$ is $\l$--independent and $\ovl\cG$ coincides with $\cG_0$ for $\varepsilon=0$, one definitely  has that, for $\varepsilon$ small, the main information on the dynamics of $\hh$ is nothing else than the one $\cG_0$, which, in the aforementioned coordinates has a very simple expression (see~\equ{G00} below).

\nl
Now, it is not simple to prove {\it directly} the convergence of the series~\equ{eq: normalization}. Therefore, we shall prove it a posteriori, under suitably more stringent assumptions, all of them verified in our application (of course, we expect that the convergence of the series~\equ{eq: normalization} holds in a more general situation). We defer this latter check to the Appendix~\ref{convergence}, being of purely technical nature.


\subsubsection{Choice of the Kepler map}\label{coordinates}

\nl
We propose a certain {\it partial Kepler map} $\textrm{\sc k}$ which includes $\L$, $\ell$ in~\equ{a} and, moreover, the functions in~\equ{trivial integrals}. The complete set of coordinates is denoted as
\beq{K}\textrm{\sc k}=(\ZZ,\GG,\Theta,\RR', \L,{\rm G}, \zz,\gamma, \vartheta, \rr', \ell, {\rm g})\eeq
where $\ZZ$ is the generalized impulse conjugated to the position coordinate  $\zz$, etc.
To define the remaining coordinates, we need the following notations. For ${\mathbf u}$, ${\mathbf v}\in\real^3$ lying in the plane orthogonal to a vector ${\mathbf w}$,  $\a_{\mathbf w}({\mathbf u},{\mathbf v})$ denotes the positively oriented angle  between ${\mathbf u}$ and ${\mathbf v}$, as seen from ${\mathbf w}$ according to the right hand rule. Define  `nodes' ${\mathbf n}_i$ as
${\mathbf n}_0:={\mathbf k}\times {\mathbf C}_{\rm t}$, ${\mathbf n}_1:={\mathbf C}_{\rm t}\times {\mathbf x}'$ and $ {\mathbf n}:={\mathbf x}'\times {\mathbf C}$ and assume that   ${\mathbf C}_{\rm t}$, ${\mathbf C}$, ${\mathbf x}'$, ${\mathbf x}$, ${\mathbf n}_0$, ${\mathbf n}_1$ and ${\mathbf n}$ do not vanish. 
Then define
\beqa{k coord}
{\rm R}':=\frac{{\mathbf y}'\cdot {\mathbf x}'}{|{\mathbf x}'|}\ ,\ 
{\rm G}:=|{\mathbf C}|\ ,\ {\rm g}:=\a_{{\mathbf C}}({\mathbf n}, {\mathbf C}\times {\mathbf P}) 
\nonumber\\
 {\rm z}:=\a_{{\mathbf k}}({\mathbf i}, {\mathbf n}_0)\ ,\ 
\gamma:=\a_{{\mathbf C}_{\rm t}}({\mathbf n}_0, {\mathbf n}_1)\ ,\ 
\vartheta:=\a_{{\mathbf x}'}({\mathbf n}_1, {\mathbf n})\ .
\eeqa

\nl
The coordinates $\textrm{\sc k}$ are canonical, since they  can be easily derived from another set of canonical coordinates\footnote{\label{changes of notations}With respect to the notations in \cite{pinzari13}, in~\equ{P coord} we have renamed $\CC_3=\ZZ$, ${\rm G}=\GG$, $\RR_1=\RR'$,   $\RR_2=\RR$,   $\zeta=\zz$, ${\rm g}=\g$, $\rr_1=\rr'$, $\rr_2=\rr$.} \beq{P coord}\textrm{\sc p}=(\ZZ, \GG, \Theta, \RR', \RR, \Phi, \zz, \gamma, \vartheta, \rr', \rr, \varphi)\eeq whose canonical character has been discussed  in \cite[Section 2]{pinzari13}, up to 
 change the quadruplet $(\RR, \Phi,\rr,\varphi)$ in $\textrm{\sc p}$ with the quadruplet $ (\L,{\rm G},\ell, {\rm g})$ in $\textrm{\sc k}$, 
 where the two quadruplets are
 related  via the classical (canonical) Delaunay map.
Note that, as a common aspect, the  coordinates
$(\ZZ,\GG, \zz)$ depend only on  ${\mathbf \CC_{\rm t}}$, while $(\Theta, \rr', \gamma)$ depend on ${\mathbf \CC_{\rm t}}$ and  ${\mathbf x}'$. Therefore all such six coordinates are first integrals to $\hh$ and ${\cal G}$. This implies that the two couples $(\ZZ,\zz)$, $(\GG,\gamma)$ and  the two coordinates ${\vartheta}$ and $\RR'$ are cyclic in $\hh_{\textrm{\sc k}}$ and ${\cal G}_{\textrm{\sc k}}$.
Then  $\hh_{\textrm{\sc k}}$ and ${\cal G}_{\textrm{\sc k}}$ have just two {\it effective} degrees of freedom, since they depend only on the six coordinates $(\L,\ell, {\rm G},{\rm g};\rr',\Theta)$, including  two only angles. Note that: (i) along the $\hh_{\textrm{\sc k}}$--motions, the angle 
 $\ell$ is a fast, while ${\rm g}$ is slow. 
Moreover, (ii) ${\cal G}_{0}$ in~\equ{G0} is $\ell$--independent, being, as said, a first integral to $\cK$ in~\equ{Kepler}. Its expression
 in terms of $\textrm{\sc k}$ is
\beqa{G00}
{\cal G}_{0}(\rr',\L,\Theta,{\rm G},{\rm g})={\rm G}^2+\mm\rr'\sqrt{1-\frac{\Theta^2}{{\rm G}^2}}\sqrt{1-\frac{{\rm G}^2}{\L^2}}\cos{\rm g}\ .
\eeqa

\subsubsection{Renormalizable integrability}\label{norm}
Let $\ovl{\textrm{\sc k}}$, ${\ovl\hh}$, ${\ovl{\cal G}}$ as in~\equ{eq: norm},~\equ{eq: normalization},~\equ{ovl G}, with  $\textrm{\sc k}$ chosen as in~\equ{K}. 
In particular,~\equ{eq: norm} becomes
\beq{eq: norm1}
\ovl{\textrm{\sc k}}=(
\ovl\Theta,\ovl\RR',\ovl\L,\ovl{\rm G},\ovl\vartheta,\ovl\rr',\ovl\ell,\ovl{\rm g}
)\to {\textrm{\sc k}}=(
\Theta,\RR',\L,{\rm G},\vartheta,\rr',\ell,{\rm g}
)
\eeq
(the four coordinates $\ZZ$, $\GG$, $\zz$, $\gamma$ have been neglected).
In terms of $\ovl{\textrm{\sc k}}$, we have
 $\ovl\hh={\ovl\hh}(\ovl\rr',\ovl\L,\ovl\Theta,\ovl{\rm G},\ovl{\rm g})$, $ \ovl{\cal G}={\ovl{\cal G}}(\ovl\rr',\ovl\L,\ovl\Theta,\ovl{\rm G},\ovl{\rm g})$. 
 Moreover, if $\cG_0$ is as in~\equ{G00},
\beq{G000}\ovl{\cal G}(\rr',\ovl\L,\ovl\Theta,\ovl{\rm G},\ovl{\rm g})=\cG_0(\ovl\rr',\ovl\L,\ovl\Theta,\ovl{\rm G},\ovl{\rm g})+{\rm O}(\varepsilon)\ .\eeq
The commutation of ${\ovl\hh}$   and ${\ovl{\cal G}}$   and  an Implicit Function Theorem argument allow to infer that  ${\ovl\hh}$ may be written as a function of ${\ovl{\cal G}}$ and the other integrals:
\beqa{ren int}{\ovl\hh}(\ovl\rr',\ovl\L,\ovl\Theta,\ovl{\rm G},\ovl{\rm g})&=&\widehat{\hh}\Big(\ovl\rr',\ovl\L,\ovl\Theta,
{\ovl{\cal G}}(\ovl\rr',\ovl\L,\ovl\Theta,\ovl{\rm G},\ovl{\rm g})\Big)\ .\eeqa
We shall use this formula in order to infer informations on the dynamics.
In view of the  application to {\sc sea}, it will be  more convenient to look at the Hamiltonian

\beqa{ovlH0}
\ovl\HH_0(\ovl\rr',\ovl\L,\ovl\Theta,\ovl{\rm G},\ovl{\rm g})&=&\widehat{\HH}_0\Big(\ovl\rr',\ovl\L,\ovl\Theta,
\ovl{\cal G}(\ovl\rr',\ovl\L,\ovl\Theta,\ovl{\rm G},\ovl{\rm g})\Big)\nonumber\\
&=&\hh_0(\ovl\rr')+\varrho\widehat\hh\Big(\ovl\rr',\ovl\L,\ovl\Theta,
\ovl{\cal G}(\ovl\rr',\ovl\L,\ovl\Theta,\ovl{\rm G},\ovl{\rm g})\Big)\nonumber\\
\eeqa
with $\hh_0$ depending only on $\ovl\rr'=\rr'$.
We then may write the solutions of $\ovl\HH_0$ in terms of the solutions of $\ovl\cG$: If $\ovl\D=\ovl\cS\times\P$
denotes the singular set $\D$ in~\equ{D} expressed in the  coordinates at left hand side of
\equ{eq: norm1}
and 
 $t\to (\ovl{\rm G}_{\ovl{\cal G}}(t),\ovl{\rm g}_{\ovl{\cal G}}(t))$
is any solution to  $(\ovl{\rm G},\ovl{\rm g})\to{\ovl{\cal G}}(\ovl\rr_0',\ovl\L_0,\ovl\Theta_0,\ovl{\rm G},\ovl{\rm g})$ with initial datum $(\ovl{\rm G}_0,\ovl{\rm g}_0)$, 
and $\ovl\rr_0',\ovl\L_0,\ovl\Theta_0$ are chosen so that 
$\big((\ovl\rr_0',\ovl\L_0,\ovl{\rm G}_0,\ovl{\rm g}_0),\ovl\Theta_0\big)\notin\ovl\cS\times\ovl\P=\ovl\D$, then one finds one solution to $\ovl\HH_0$ by letting
\beqa{solutions0}
&&\arr{\ovl\rr'(t)=\ovl\rr'_0\\\\
 \ovl\L(t)=\ovl\L_0\\\\
 \ovl\Theta(t)=\ovl\Theta_0
 }\ \ \arr{\ovl{\rm G}(t)=\ovl{\rm G}_{\ovl{\cal G}}(\varrho\widehat\hh_{\ovl\cG} t)\\\\
 \ovl{\rm g}(t)=\ovl{\rm g}_{\ovl{\cal G}}(\varrho\widehat\hh_{\ovl\cG} t)}\nonumber\\
&& \arr{\ovl\RR'(t)=\ovl\RR_0'-(\hh_0'+\varrho\widehat\hh_{\ovl\rr'}+\varrho\widehat\hh_{\ovl\cG}\ovl\cG_{\rr'})(t-t_0)\\\\
\ovl\ell(t)=\ovl\ell_0+\varrho(\widehat\hh_{\ovl\L}+\widehat\hh_{\ovl\cG}\ovl\cG_\L)(t-t_0)\\\\
\ovl\vartheta(t)=\ovl\vartheta_0+\varrho(\widehat\hh_{\ovl\Theta}+\widehat\hh_{\ovl\cG}\ovl\cG_\Theta)(t-t_0)
}
 \eeqa
 where the derivatives are evaluated at $\big(\ovl\rr_0'$, $\ovl\L_0$, $\ovl\Theta_0$,
$\ovl{\cal G}(\ovl\rr_0'$, $\ovl\L_0$, $\ovl\Theta_0$, $\ovl{\rm G}_0$, $\ovl{\rm g}_0)\big)$.

\subsubsection{First order  solutions in $\varepsilon$}
It is possible to prove (see Appendix~\ref{lambda independence}) the following `local' version of~\equ{ren int}
\beq{U}\ovl\UU(\ovl\rr',\ovl\L,\ovl\Theta,\ovl{\rm G},\ovl{\rm g})=\widehat\UU\Big(\ovl\rr',\ovl\L,\ovl\Theta,{{\cal G}_{\textrm{0}}}(\ovl\rr',\ovl\L,\ovl\Theta,\ovl{\rm G},\ovl{\rm g})\Big)\eeq
where ${{\cal G}_{0}}$ is as in~\equ{G00}. 
This formula  turns to be more  useful than~\equ{ren int} from the practical point of view  because now
the function $\widehat\UU$ may be {\it explicitly} written (see Appendix~\ref{explicit}). 
We then have a
 concrete first order approximation of~\equ{solutions0}:
  \beqa{solutions}
 &&\left.
 \begin{array}{lll}\ovl\rr'(t)=\ovl\rr'_0\\\\
 \ovl\L(t)=\ovl\L_0\\\\
 \ovl\Theta(t)=\ovl\Theta_0
  \end{array}
 \right\}\ \left.
 \begin{array}{llll}\ovl{\rm G}(t)=\ovl{\rm G}_0(\varrho\varepsilon\widehat\UU_{{\cG_0}} t)\\\\
 \ovl{\rm g}(t)=\ovl{\rm g}_0(\varrho\varepsilon\widehat\UU_{{\cG_0}} t)
 \end{array}
 \right\}+{\rm O}(\varepsilon; \ovl\rr',\ovl\L,\ovl\Theta, \varrho\varepsilon t)\nonumber\\
&&\left. \begin{array}{llll}\ovl\RR'(t)=\ovl\RR_0'-(\hh_0'+\varrho\varepsilon\widehat\UU_{\rr'}\\
\qquad\ \ +\varrho\varepsilon\widehat\UU_{{\cG_0}}{\cG_0}_{\ovl\rr'})(t-t_0)\\\\
\ovl\ell(t)=\ovl\ell_0+\varrho\varepsilon(\widehat\UU_{\L}\\
\qquad\ \ +\widehat\UU_{{\cG_0}}{\cG_0}_{\ovl\L})(t-t_0)\\\\
\ovl\vartheta(t)=\ovl\vartheta_0+\varrho\varepsilon(\widehat\UU_{\Theta}\\
\qquad\ \ +\widehat\UU_{{\cG_0}}{\cG_0}_{\ovl\Theta})(t-t_0)
 \end{array}
 \right\}+{\rm O}(\varepsilon; \ovl\rr',\ovl\L,\ovl\Theta, \varrho\varepsilon t)\eeqa
where $t\to(\ovl{\rm G}_0(t),\ovl{\rm g}_0(t))$ is any solution to $(\ovl{\rm G},\ovl{\rm g})\to\cG_0(\ovl\rr_0',\ovl\L_0,\ovl\Theta_0,\ovl{\rm G},\ovl{\rm g})$.\\
According to~\equ{solutions}, the problem is reduced to study the dynamics of $\cG_0$ in~\equ{G00}, for $\textrm{\sc k}\in \cD$.
The study of the phase portrait of $\cG_0$ reduces to study an equation with degree three in ${\rm G}^2$ (since equation~\equ{G00} can be written in this way). For such an equation, as well known, one could use Cardano formulae. However,
 the purpose of this paper is to highlight the utility of the method, rather than to push the analysis to its maximum generality.
 Therefore, we
 shall simplify even more the discussion restricting to $\P=\{\Theta=0\}$, in which case the equation to be solved is of degree two, completely explicit.

\begin{remark}\label{rem}
Equations~\equ{solutions}  imply that on the 
manifolds defined by values of $\ovl\rr'_0$, $\ovl\L_0$, $\ovl\Theta_0$, $\ovl\cG$ where $\widehat\hh_{\ovl\cG}$ vanishes,  the coordinates  $(\ovl{\rm G}$, $\ovl{\rm g})$ do not move. In view of~\equ{eq: normalization} and~\equ{U} above, and the Implicit Function Theorem, such zeroes may be  approximated by the corresponding zeroes of $\widehat\UU_{\cG_0}$.
From the analysis of an $\frac{a}{\rr'}$--expansion  \beq{hatU}\widehat\UU=\frac{1}{a}\Big(1-\frac{{\rr'}^2}{4a^2}\frac{\L^3(3\Theta^2-{\cal G}_0^2)}{{\cal G}_0^5}+\cdots\Big)\eeq
 one can see that, as a matter of fact, $\widehat\UU_{\cG_0}$ does vanish for
\beq{fixed points}\cG_0\sim \pm\sqrt5\Theta\ .\eeq
Thus, when $\rr'\ll a$,~\equ{fixed points} is a first approximation of manifolds where $(\ovl{\rm G},\ovl{\rm g})$ stay fixed for $\ovl\hh$.
\end{remark}
 
 \subsubsection{Equilibria and  phase portrait in the planar case}\label{G0levels}

\begin{figure}
 \includegraphics[height=4.0cm, width=8.0cm]{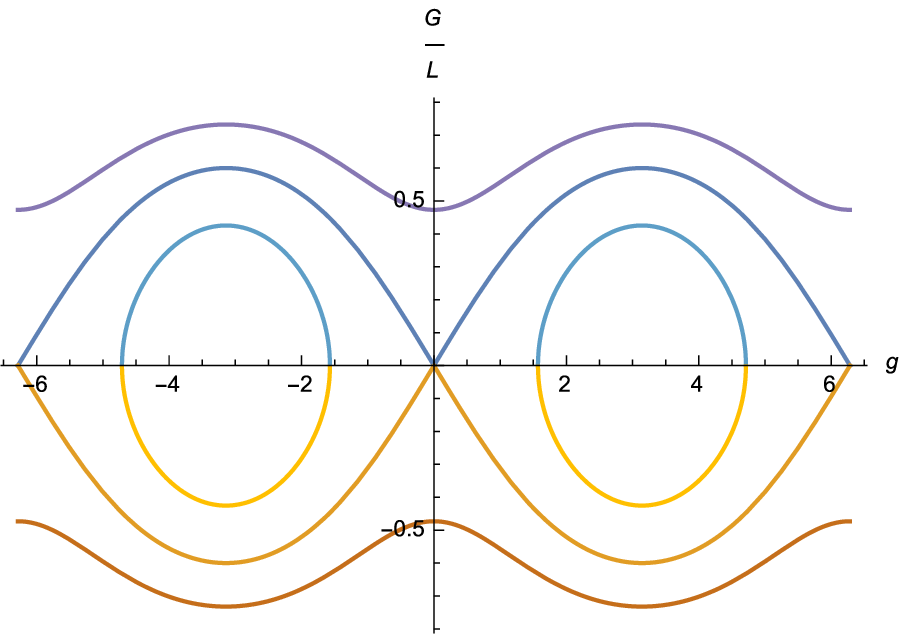}%
 \caption{The phase portrait of $\widehat\cG_0$, for $0<\d<1$ and $-\d\le\widehat\cG_0<1$ ({\sc mathematica}).
  }\label{fig: stable}
 \end{figure}

It is convenient to divide equation~\equ{G00} by $\L^2$ and let $\widehat\cG_0:=\frac{\cG_0}{\L^2}$, $\d:=\frac{\rr'}{a}$. 
Then equation~\equ{G00} with  $\Theta=0$ 
becomes
\beq{G0pl}\widehat\cG_0=\frac{{\rm G}^2}{\L^2}+\d \sqrt{1-\frac{{\rm G}^2}{\L^2}}\cos{\rm g}\ . \eeq
 Let us first discuss the equilibria of $\widehat\cG_0$.  This function is even and regular around  $({\rm g},{\rm G})=(\p,0)$ and
$({\rm g},{\rm G})=(0,0)$ (note that for $\Theta\ne0$, ${\rm G}=0$ would be a singularity), which are equilibria. However, the character of such equilibria  is different accordingly to wether
 $\d\in (2,+\infty)$ or $\d\in (0,2)$: in the former case, $(0,0)$ and $(\p,0)$ are both  stable; in the latter case, $(0,0)$ is unstable, while $(\p,0)$ is stable.  Correspondingly to such extrema,  $\widehat\cG_0$ takes the values
 $\widehat\cG_{0\rm sad}=\d$ and $\widehat\cG_{0\rm min}=-\d$, respectively.
 To such equilibria, one should add, for $\d<2$,  also the point $(0, \sqrt{1-\frac{\d^2}{4}})$, which is stable, and where $\widehat\cG_0$ takes the maximum possible vale $\widehat\cG_{0\rm max}=1+\frac{\d^2}{4}$.
We now study the phase portrait of $\widehat\cG_0$ in~\equ{G0pl}. Still, we shall make a further simplification. We take $0<\d<1$ and $-\d\le \widehat\cG_0\le1$, leaving the remaining analysis to the interested reader. The level $\widehat\cG_0=1$ splits as $$\{\widehat\cG_0=1\}=\{{\rm G}=\L\}\bigcup\{{\rm G}=\L|\sin{\rm g}|\}\ .$$
For $-\d\le\widehat\cG_0<1$, we rewrite equation 
\equ{G0pl}
in terms of $w:=\sqrt{1-\frac{{\rm G}^2}{\L^2}}$, with $w\in (0,1)$. We obtain
\beq{equal}w^2-\d w\cos{\rm g}-1+\widehat\cG_0=0\ .\eeq
Solving for $w$, 
\beq{wpm}w_\pm=\frac{\d\cos{\rm g}\pm\sqrt{\d^2\cos^2{\rm g}+4-4\widehat\cG_0}}{2}\eeq
we see that $w_-$ is negative for all ${\rm g}$, and hence is to be disregarded, while $w_+$ is positive for all ${\rm g}$. Moreover, $w_+$ does not exceed 1 as soon as
 \beq{ineq}\sqrt{\d^2\cos^2{\rm g}+4-4\widehat\cG_0}<2-\d\cos{\rm g}\ .\eeq
Solving for $\cos{\rm g}$, we obtain
 \beq{ineq1}\cos{\rm g}< \min\Big\{\frac{2}{\d},\  \frac{\widehat\cG_0}{\d}\Big\}=\frac{\widehat\cG_0}{\d}\ .\eeq
From this inequality we see that if $\widehat\cG_0>\d$, we have a rotational motion of ${\rm g}$, while if $-\d\le \widehat\cG_0<\d$ the motion is librational around $(\p,0)$. The situation is depicted in Figure~\ref{fig: stable}.
 
 \subsubsection{The separatrix}
The  separatrix  in Figure~\ref{fig: stable}, namely the $\widehat\cG_0$--level though $(0,0)$, having equation 
   $\widehat\cG_0=\d$,
    is precisely the set $\cS$ in~\equ{S}. Therefore, the motions along it, strictly speaking, loose their meaning if replaced within the formulae~\equ{solutions}.
Note  that the homoclinic solution  along $\cS$ can however be  easily computed: 
\beq{Gammagamma}
\arr{\dst
\ovl{\rm G}_0(t)=\frac{\s\ovl\L}{\cosh\s\ovl\L(t-t_0)}\\\\
\dst\ovl{\rm g}_0(t)=\pm\cos^{-1}\frac{1-\frac{\b^2}{\cosh^2\s\ovl\L(t-t_0)}}{\sqrt{1-\frac{\s^2}{\cosh^2\s\ovl\L(t-t_0)}}}
}\eeq
with
$\s^2:=\d(2-\d)$, $\b^2:=2-\d$.

\subsubsection{Action--Angle coordinates}

The previous analysis shows that, by Liouville--Arnold  \cite{arnold64}, for all $0<\d<1$ in any of the two regions connected  region  defined by inequalities $-\d\le \widehat\cG_0<\d$, $\d<\widehat\cG_0<1$,
it is possible to find a (different) canonical change of coordinates
\beq{eq: norm2}
{\textrm{\sc a}}=(
\widehat\RR',\widehat\L,\widehat A,\widehat\rr',\widehat\ell, \widehat\a)\to \ovl{\textrm{\sc k}}=(
\ovl\RR',\ovl\L,\ovl{\rm G},\ovl\rr',\ovl\ell,\ovl{\rm g}
)
\eeq
(preserving $\widehat\L=\ovl\L$, $\widehat\rr'=\ovl\rr'$) with $\widehat\a\in \torus$, such that, in their terms,  $\cG_0$ becomes a function of only $(\widehat\rr',\widehat\L,\widehat A)$.
By usual  integrability arguments, up to $\varepsilon$--variations of the coordinates, one has also $\cG$ is so and hence also, by~\equ{ovlH0}, $\HH_0$, with all these functions
being well defined (collision--free).

\subsection{Application to the Sun--Earth--Asteroid  system}
Let us  consider  the problem of three gravitational masses, $1$, $\varepsilon$,  $\m$, with $1\gg \varepsilon\gg \m$.
After the reduction of translation invariance according to the heliocentric method, the three--body problem in $\real^3$ is governed by the six--degrees--of--freedom Hamiltonian

\beqa{3B}\ovl\HH(\ovl {\mathbf y}',\ovl {\mathbf y},\ovl {\mathbf x}',\ovl {\mathbf x})&=&\frac{|\ovl {\mathbf y}'|^2}{2\varepsilon {\mm}'}-\frac{\varepsilon}{|\ovl {\mathbf x}'|}+\frac{|\ovl {\mathbf y}|^2}{2\m {\mm}}-\frac{\m}{|\ovl {\mathbf x}|}\nonumber\\
&&-\frac{\m\varepsilon  }{|\ovl {\mathbf x}'-\ovl {\mathbf x}|}+{\ovl {\mathbf y}'\cdot \ovl {\mathbf y}}{}
\eeqa
where $\mm':=(1+\varepsilon)^{-1}$, 
$\mm:=(1+\m)^{-1}$
;
 $\ovl {\mathbf y}$, $\ovl {\mathbf y}'$, $\ovl {\mathbf x}$, $\ovl {\mathbf x}'$ $\in$ $\real^3$ are impulse--position coordinates.

\nl
In order to eliminate small numbers from  denominators, one rescales time, Hamiltonian and coordinates, via 
$
\HH( {\mathbf y}', {\mathbf y}, {\mathbf x}', {\mathbf x}):=\varepsilon^{-1}\ovl\HH(\m  {\mathbf y}',\m  {\mathbf y}, {\mathbf x}', {\mathbf x})$
so that $\HH$ describes the evolution of
${\mathbf y}'$, ${\mathbf y}$, ${\mathbf x}'$, ${\mathbf x}$
during the time $\varrho^{-1} t$, where $\varrho:=\m/\varepsilon$.

\nl
We obtain
\beqa{h3B}\HH={\rm h}_0
+\varrho{\rm h}_{1}+ \varrho^2{f}=\HH_0+\varrho^2 f\eeqa
where $\HH_0:={\rm h}_0
+\varrho{\rm h}$, ${f}:=\frac{|{\mathbf y}'|^2}{2{\mm}'}+\varepsilon{\mathbf y}'\cdot {\mathbf y}$ with $\hh_0:=-\frac{1}{|{\mathbf x}'|}$. Let us write such functions  in terms of $\textrm{\sc k}$, without changing them the names. 
Since $f$   still possesses ${\mathbf C}_{\rm t}$ as a first integral, but no longer ${\mathbf x}'$, we have that  $f$ is a function of 
 $\RR'$,  $\L$, ${\rm G}$,  $\Theta$, $\rr'$, $\ell$, ${\rm g}$, $\vartheta$ and, moreover, it depends {\it parametrically}, on $\GG$. 
The manifolds
$\P_{\downarrow}:=\big\{(\Theta,\vartheta)=(0,0)\big\}$ and $\P_{\uparrow}:=\big\{(\Theta,\vartheta)=(0,\p)\big\}$
are invariant to the $\HH$--flow. 
Motions on $\P_{\uparrow}$, $\P_{\downarrow}$ correspond  to have, at all times, the orbits of Earth and the Sun on the same instantaneous plane, with a suitable choice of the mutual  inclination ($0$ or $\p$) of ${\mathbf C}_{\rm t}$, ${\mathbf C}$.
We  turn to the coordinates $\ovl{\textrm{\sc k}}$ in~\equ{eq: norm1}. Then $\HH_{0}$ is carried to $\ovl\HH_0$ in~\equ{ovlH0}, while $f$ to a suitable $\ovl f$.
The dynamics of $\ovl\HH_0$ on  $\P$ (and hence, in particular, on its sub--manifolds $\P_\uparrow$, $\P_\downarrow$) has been discussed in the previous paragraph. We consider the motions corresponding to the cases (a), (b) or (c), and aim to extend (many of) them to {\sc sea}. We put the system in the coordinates $\textrm{\sc a}$ in~\equ{eq: norm2}, letting  $\widehat\HH({\widehat\RR}',\widehat\L,\widehat A,\widehat\rr',\widehat\ell,\widehat\a)=\widehat\HH_0(\widehat\L,\widehat A,\widehat\rr')+\varrho^2\widehat f({\widehat\RR}',\widehat\L,\widehat A,\widehat\rr',\widehat\ell,\widehat\a) 
$ the relative Hamiltonian. Via  normal form theory  (see Appendix~\ref{normal form}), we conjugate $\widehat\HH$ to (omitting to write the dependence on $\GG$) \beq{tildeH} \widetilde\HH=\widehat\HH_0(\widetilde\L,\widetilde A,\widetilde\rr')+\varrho^2\widetilde\HH_1({\widetilde\RR}',\widetilde\L,\widetilde A,\widetilde\rr')+\varrho^2\widetilde f({\widetilde\RR}',\widetilde\L,\widetilde A,\widetilde\rr',\widetilde\ell,\widetilde a)
\eeq
where  $\widetilde f$ is of order $2^{-c/\varrho}$, while, up to higher orders,
 \beq{2body}\widetilde\HH_1\sim\frac{({\widetilde\RR}')^2}{2\mm'}+\frac{\Phi(\widetilde\L,\widetilde A,\widetilde \rr')^2}{2\mm'{\rr'}^2}\eeq
with $\Phi(\widetilde\L,\widetilde A,\widetilde \rr')$ a suitable regular function.
At this point, one  integrates the term $\widehat\HH_0+\varrho^2\widetilde\HH_1$ 
with respect to $({\widetilde\RR}',{\widetilde\rr}')$. Since, by construction, $\widehat\HH_0$ is regular, and because of~\equ{2body}, the integration of such term is analogous to a two--body system, and one finds, for low energies, an action--angle couple $(\L',\ell')$ such that $\widehat\HH_0+\varrho^2\widetilde\HH_1$, after the integration, would depend on $(\widetilde\L, \widetilde\L',\widetilde A)$ only. An application of {\sc kam} theory allows to infer  the existence of quasi--periodic motions with three frequencies, with a residual set having an exponentially small density.

\newpage\section{Conclusions and perspectives}\label{perspectives}
We proposed a new analysis of the two-centre Hamiltonian $\hh$ in~\equ{2C} based on the Euler integral $\cG$ in~\equ{cal G}. We introduced  a `ad--hoc' system of canonical coordinates~\equ{K} which includes all its first integrals but $\cG$.
Accordingly, we wrote $\hh$  as an effective, properly degenerate, two--degrees of freedom system, with a fast angle $\ell$ and a slow one ${\rm g}$. Eliminating (in the regime where the two attracting centres have much different masses ratio $\varepsilon$) the fast angle $\ell$ via perturbative methods, we obtained a new one degrees of freedom Hamiltonian $\ovl\hh$ in~\equ{eq: normalization} and proved that the dynamics of $\ovl\hh$ is completely determined by the one of $\cG$, written in the new coordinates, the function $\ovl\cG$ in~\equ{ovl G}. The result carries an important consequence: at the lowest  order in $\varepsilon$, the common motions of $\ovl\hh$ and $\ovl\cG$ are determined by the simple one--dimensional Hamiltonian $\cG_0$ in~\equ{G00}. In 
 the case of the planar problem,
the phase portrait of $\cG_0$ is  explicitly, rigorously computable (Section~\ref{G0levels} and Figure~\ref{fig: stable}). It shows that Liouville--Arnold
 action--angle coordinates $(A,a)\iff ({\rm G},{\rm g})$  do exist for all $\rr'<a$ and $-\mm\rr'\le \cG_0<\L^2$, apart for a zero measure set $\cS$, where collisions are possible. We applied the result 
to the {\sc sea} system, regarding it as a perturbation of 2{\sc cp}.
 A suitable normal form theory and {\sc kam} theory allows to infer the existence of quasi--periodic motions with three frequencies, with a residual set having an exponentially small density. As a byproduct of the proof, the risk of Earth--Asteroid collision may be excluded in all cases where $\cG_0$ is sufficiently far away  from $\mm\rr'$, where $\mm$ is the Earth mass, in suitable units, and $\rr'$ is distance from the sun. 

\nl
We conjecture that closely to the collision sets $\cS$'s, {\sc sea} possesses chains of {\it transition tori} in the sense of   \cite{arnold64}.

\section*{Acknowledgements}

Figure~\ref{fig: stable} was produced with  {\sc mathematica}. The author thanks M. Guzzo for his encouragement.



\appendix

\newpage\section{Technical details to Section 2}

\subsection{ The  Euler integral}\label{App: 2C}
The formulae in~\equ{cal G}--\equ{G0} 
are not standard in the literature.  For sake of completeness, and for the reader's facility, we report their derivation  here.

\nl
In Section~\ref{proof of G1}, we check that, writing the two--centre Hamiltonian in the more usual `symmetric' form
\beq{2Cnew}\hh_{\rm sim}=\frac{|{\mathbf y}|^2}{2}-\frac{m_+}{|{\mathbf x}+{\mathbf x}_0|}-\frac{m_-}{|{\mathbf x}-{\mathbf x}_0|}\eeq
then the Euler integral to $\hh_{\rm sim}$ is given by
\beq{G1}\cG_{\rm sim}=|{\mathbf x}\times {\mathbf y}|^2+({\mathbf x}_0\cdot {\mathbf y})^2+2 {\mathbf x}\cdot {\mathbf x}_0\big(\frac{m_+}{|{\mathbf x}+{\mathbf x}_0|}-\frac{m_-}{|{\mathbf x}-{\mathbf x}_0|}\big)\ .\eeq
In  Section~\ref{derivation}, we shall check that, when $\hh$ is written in the form~\equ{2C}, then $\cG_{\rm sim}$ reduces to $\cG$ in~\equ{cal G}--\equ{G0}.

\nl
Observe, incidentally, that, in the symmetric case, when the two stars merge, e.g., ${\mathbf x}_0=0$, ${\cal G}_{\rm sim}$ reduces to $|{\mathbf C}|^2$.

\nl
\subsubsection{Derivation of~\equ{G1}}\label{proof of G1}
For part of the proof, we use the canonical coordinates $\textrm{\sc p}$ in~\equ{P coord} (for uniformity of notations, we shall  replace the symbols ${\mathbf x}'$, ${\mathbf y}'$, $\rr'$, $\RR'$ in~\equ{k coord} and~\equ{P coord} with ${\mathbf x}_0$, ${\mathbf y}_0$, $\rr_0$, $\RR_0$, respectively).
As said, the $\textrm{\sc p}$'s  have in common with $\textrm{\sc k}$'s in~\equ{K} almost all the coordinates, apart for the two quadruplets $(\RR$, $\Phi$, $\rr$, $\f)$ (for the $\textrm{\sc p}$'s) and $(\L$, ${\rm G}$, $\ell$, ${\rm g})$ (for the $\textrm{\sc k}$'s).
The definition of the former is\footnote{The angle $\f$ in~\equ{quadruplet} corresponds to $\f+\frac{\p}{2}$ of \cite[equation (2.10)]{pinzari13}.} (within the same notations as in Section~\ref{coordinates})
\beq{quadruplet}\RR=\frac{{\mathbf y}\cdot {\mathbf x}}{|{\mathbf x}|}\ ,\ \Phi=|{\mathbf C}|\ ,\ \rr=|{\mathbf x}|\ ,\ \f=\a_{\mathbf C}({\mathbf n},{\mathbf k}\times {\mathbf x})\ .\eeq
 In terms of  $\textrm{\sc p}$, the scalar product ${\mathbf x}_0\cdot {\mathbf x}$
takes the form
$${ {\mathbf x}_0\cdot {\mathbf x}=-\rr_0 \rr \sqrt{1-\frac{\Theta^2}{\Phi^2}}\cos\varphi}$$
and so $\hh_{\rm sim}$ in~\equ{2C} becomes
\beq{H2CRPhi}\hh_{\rm sim}=\frac{{\rm R}^2}{2}+\frac{\Phi^2}{2 \rr^2}-\frac{m_+}{\rr_+}-\frac{m_-}{\rr_-}\eeq
where
$${\rr_\pm^2:=\rr_0^2\mp2 \rr_0 \rr \sqrt{1-\frac{\Theta^2}{\Phi^2}}\cos\varphi +\rr^2\ .}$$
$\hh_{\rm sim}$ has now two degrees of freedom, exactly as in the classical discussion, which goes along the  `elliptic coordinates' 
$$\l=\frac{\rr_++\rr_-}{2 \rr_0}\qquad \m=\frac{\rr_+-\rr_-}{2\rr_0}\ .$$
 We then define a  change of canonical coordinates $(\RR,\Phi,\rr,\f)\to (p_\l,p_\m,\l,\m)$ where $\l$, $\m$ are as above, while  their conjugated momenta $p_\l$, $p_\m$ are found taking the inverse of
\beq{rpm}\rr_+=\rr_0(\l+\m)\qquad \rr_-=\rr_0(\l-\m)\eeq
and than squaring and summing, or subtracting. This gives
\beq{r}{\rr=\rr_0\sqrt{\l^2+\m^2-1}\qquad \varphi=\cos^{-1}\Big(-\frac{\l\m}{\sqrt{\l^2+\m^2-1}\sqrt{1-\frac{\Theta^2}{\Phi^2}}}\Big)}\eeq
Then one considers the generating function
\beqano
S(\Phi, \Theta, {\rm R}_0,{\rm R}, \l,\m, \hat \rr_0, \hat\vartheta)&=&\Theta\hat\vartheta+{\rm R}_0 \hat \rr_0+{\rm R}\hat \rr_0\sqrt{\l^2+\m^2-1} \nonumber\\
&&+\int^\Phi  \cos^{-1}\nonumber\\
&&\Big(-\frac{\l\m}{\sqrt{\l^2+\m^2-1}\sqrt{1-\frac{\Theta^2}{\Phi^2}}}\Big)d\Phi\ .
\eeqano
The transformation generated by $S$ leaves the coordinates
$\Theta$, $\rr_0$
unvaried (therefore, we shall not change their names), while shifts in an inessential way (since they do not appear into $\hh_{\rm sim}$) the coordinates $\vartheta$, ${\rm R}_0$. Taking the derivatives with respect to $\l$, $\m$, one finds
$$\arr{
\dst p_\l=\frac{\rr_0\l {\rm R}}{\sqrt{\l^2+\m^2-1}}-\frac{\m\sqrt{(1-\m^2)(\l^2-1)\Phi^2-(\l^2+\m^2-1)	\Theta^2}}{(\l^2+\m^2-1)(\l^2-1)}\\\\
\dst p_\m=\frac{\rr_0\m {\rm R}}{\sqrt{\l^2+\m^2-1}}+\frac{\l \sqrt{(1-\m^2)(\l^2-1)\Phi^2-(\l^2+\m^2-1)	\Theta^2}}{(\l^2+\m^2-1)(1-\m^2)}
}
$$
whence, taking the inverse with respect to ${\rm R}$, $\Phi$
$$
\arr{\dst{\rm R}=\frac{\l(\l^2-1)p_\l+\m(1-\m^2)p_\m}{\rr_0(\l^2-\m^2)\sqrt{\l^2+\m^2-1}}\\\\
\dst \Phi^2=\frac{(\l p_\m-\m p_\l)^2(\l^2-1)(1-\m^2)}{(\l^2-\m^2)}+\frac{\l^2+\m^2-1}{(1-\m^2)(\l^2-1)}\Theta^2}
$$
Replacing these expressions and the one for  $\rr_+$, $\rr_-$, $\rr$ in~\equ{rpm},~\equ{r} into the Hamiltonian $\hh_{\rm sim}$ in~\equ{H2CRPhi}, one finds the classical expression
\beqa{H2C}\hh_{\rm sim}&=&\frac{p^2_\l(\l^2-1)}{2 \rr_0^2(\l^2-\m^2)}+\frac{p^2_\m(1-\m^2)}{2 \rr_0^2(\l^2-\m^2)}\nonumber\\
&&+\frac{\Theta^2}{2 \rr_0^2(\l^2-\m^2)}\big(\frac{1}{1-\m^2}+\frac{1}{\l^2-1}\big)\nonumber\\
&&-\frac{(m_++m_-)\l-(m_+-m_-)\m}{\rr_0^2(\l^2-\m^2)}\ .
\eeqa
Then one sees that Hamilton--Jacobi equation
$$\hh_{\rm sim}-{\rm E}=0$$
splits as
\beq{split}{\cal F}^{(\mu)}(p_\m,\m,\Theta, {\rm E}, \rr_0)-{\cal F}^{(\l)}(p_\l,\l,\Theta, {\rm E}, \rr_0)=0\eeq
where
$${\cal F}^{(\m)}=p^2_\m(1-\m^2)+\frac{\Theta^2}{1-\m^2}+2(m_+-m_-)\m+2 \rr_0^2\m^2{\rm E}$$
$${\cal F}^{(\l)}=-p^2_\l(\l^2-1)-\frac{\Theta^2}{\l^2-1}+2(m_++m_-)\l+2 \rr_0^2\l^2{\rm E}\ .$$

\nl
Equation ~\equ{split} implies then that ${\cal F}^{(\mu)}(p_\m,\m,\Theta, {\rm E}, \rr_0)=\cG_{\rm sim}^{(\m)}(\Theta, {\rm E}, \rr_0)$ is actually independent of $(p_\m,\m)$;  ${\cal F}^{(\l)}(p_\m,\m,\Theta, {\rm E}, \rr_0)=\cG_{\rm sim}^{(\l)}(\Theta, {\rm E}, \rr_0)$ is actually independent of $(p_\l,\l)$, and, a fortiori, since the partial derivatives of $\cG_{\rm sim}^{(\m)}$, $\cG_{\rm sim}^{(\l)}$ depend explicitly on $\m$, $\l$, there must exists a $\cG_{\rm sim}\in \real$ such that
$${\cal F}^{(\m)}={\cal F}^{(\l)}=\cG_{\rm sim}\ .$$
Therefore,
\beqano
\cG_{\rm sim}&=&\frac{1}{2}({\cal F}^{(\m)}+{\cal F}^{(\l)})\nonumber\\
&=&\frac{p^2_\m}{2}(1-\m^2)-\frac{p^2_\l}{2}(\l^2-1)+\frac{\Theta^2}{2}\big(\frac{1}{1-\m^2}-\frac{1}{\l^2-1}\big)\nonumber\\
&&+m_+(\l+\m)+m_-(\l-\m)\nonumber\\
&&+2\rr_0^2(\l^2+\m^2){\rm E}\ .
\eeqano
After some elementary computations, one finds the expression of  $\cG_{\rm sim}$ in terms of the coordinates $\textrm{\sc p}$  is
\beqano
\cG_{\rm sim}&=&\Phi^2+\rr_0^2(1-\frac{\Theta^2}{\Phi^2})(-{\rm R}\cos\varphi+\frac{\Phi}{r}\sin\varphi)^2\nonumber\\
&&- 2 \rr \rr_0\cos\varphi\sqrt{1-\frac{\Theta^2}{\Phi^2}}\big(\frac{m_+}{\rr_+}-\frac{m_-}{\rr_-}\big)
\eeqano
While, in terms of the coordinates $({\mathbf y}_0, {\mathbf x}_0)$, $({\mathbf y}, {\mathbf x})$, $\cG_{\rm sim}$ has the expression  in~\equ{G1}.

\subsubsection{Derivation of~\equ{cal G}--\equ{G0}}\label{derivation}

Let $\hh$ be as in~\equ{2C}.
We  preliminarily rescale   $\hh$, letting
\beqa{rescale}
\widehat\hh(\widehat {\mathbf y}, \widehat {\mathbf x}, \widehat {\mathbf x}')&=& {\rm m}^{-1}\hh({\rm m}\widehat {\mathbf y}, \widehat {\mathbf x}, \widehat {\mathbf x}')\nonumber\\
&=&\frac{| \widehat  {\mathbf y}|^2}{2}-\frac{\mm^{-1}}{|\widehat {\mathbf x}|}-\frac{\varepsilon  \mm^{-1}}{|\widehat {\mathbf x}'- \widehat {\mathbf x}|}\ .\eeqa
Letting further
\beq{helio}
 \widehat {\mathbf y}'=\frac{1}{2}({\mathbf y}_0-\ovl{\mathbf y})\qquad\widehat {\mathbf y}=\ovl{\mathbf y}\qquad  \widehat {\mathbf x}'=2{\mathbf x}_0\qquad\widehat {\mathbf x}={\mathbf x}_0+\ovl{\mathbf x}
\eeq
we approach the   Hamiltonian $\hh_{\rm sim}$ in~\equ{2Cnew}, with masses
$$m_-=\varepsilon\mm^{-1}\qquad m_+=\mm^{-1}$$
and $({\mathbf y},{\mathbf x})$ replaced by $(\ovl{\mathbf y},\ovl{\mathbf x})$. But $\hh_{\rm sim}$ admits the integral ${\cal G}_{\rm sim}$ in~\equ{G1}, and hence, applying the inverse transformations of~\equ{helio}
and~\equ{rescale}
 we find that  $\hh$ 
has the first integral
\beqano
\frac{\widehat\cG}{\mm}&:=&\frac{1}{{\rm m}}|({\mathbf x}-\frac{{\mathbf x}'}{2})\times {\mathbf y}|^2+\frac{1}{4{\rm m}}({\mathbf x}'\cdot {\mathbf y})^2\nonumber\\
&+&{\mathbf x}'\cdot({\mathbf x}-\frac{{\mathbf x}'}{2})\big(\frac{1}{|{\mathbf x}|}-\frac{\varepsilon }{|{\mathbf x}'-{\mathbf x}|}\big)
\eeqano
After multiplying  by $\mm$, we rewrite this integral as
\beqa{cal N}
\widehat\cG=\cG_0+\varepsilon\cG_1+\mm\frac{|{\mathbf x}'|^2}{2}{\hh}\eeqa
where
\beqano
\cG_0:=|{\mathbf C}|^2-{\mathbf x}'\cdot {\mathbf L}\qquad 
\cG_1:= \mm\frac{({\mathbf x}'-{\mathbf x})\cdot {\mathbf x}'}{|{\mathbf x}'-{\mathbf x}|}\eeqano
with
\beqno{\mathbf C}:={\mathbf x}\times {\mathbf y}\qquad {\mathbf L}={\mathbf y}\times{\mathbf C}-\mm\frac{{\mathbf x}}{|{\mathbf x}|}\ .\eeqno

\nl
Since the last term in~\equ{cal N} is  itself an integral for $\hh$, we can neglect it and conclude that the function
\beqano
\cG:=\cG_0+\varepsilon\cG_1\eeqano
is an integral to $\hh$. We   recognize that ${\mathbf C}$, ${\mathbf L}$ are the angular momentum and  the eccentricity vector associated to 
${\rm h}$, respectively.  This is exactly  what we had to check, after recalling that ${\mathbf L}$ is related (in our units)  to $\ee$ and ${\mathbf P}$ via $ {\mathbf L}=\mm {\rm e} {\mathbf P}$.

\subsection{On the convergence of the series~\equ{eq: normalization}}\label{convergence}

Since two different partial Kepler maps $\textrm{\sc k}=\big((\L,u,v),\l\big)$, $\textrm{\sc k}'=\big((\L,u',v'),\l'\big)$ are linked by a relation of the form
\beq{relation}
\L=\L ,(u,v)=F(\L, u',v') ,\ \l=\l'+\psi(\L,u',v')
\eeq
the character of the series~\equ{eq: normalization} does not depend on the choice of $\textrm{\sc k}$. Therefore, we choose $\textrm{\sc k}$ as in Section~\ref{coordinates}.
We prove that the series~\equ{eq: normalization} converges in the domain defined by the following inequalities
\beq{new domain}\ovl\cD:\ \Theta=0\ ,\ 0<\d<1\ ,\ \cG_0\in [-\d,\d)\cup(\d,1)\ ,\eeq
which is enough for our purposes.
With this choice,  $\hh_{\textrm{\sc k}}$ and $\cG_{\textrm{\sc k}}$ depend, as already remarked, only on the two angles $\l$ and ${\rm g}$. Let us discuss the question using Liouville--Arnold theorem.
Regarding $\hh_{\textrm{\sc k}}$ and $\cG_{\textrm{\sc k}}$ as functions of $(\L,{\rm G}, \ell, {\rm g})$, we look at level sets
$$\cM_{\hh, \ovl\cG,\varepsilon}=\Big\{(\L,{\rm G}, \ell, {\rm g},\varepsilon):\ \hh_{\textrm{\sc k}}(\L,{\rm G}, \ell, {\rm g},\varepsilon)=\hh,\ \cG_{\textrm{\sc k}}(\L,{\rm G}, \ell, {\rm g})=\ovl\cG\Big\}\ .$$
For $\varepsilon=0$, $\hh_{\textrm{\sc k}}$ reduces to $\hh_\cK$ in~\equ{Kepler}, while $\cG$ reduces to $\cG_0$ in~\equ{G0pl}. Therefore, $\cM_{\hh, \ovl\cG,0}$ is the product $\big\{\hh_\cK(\L)=\hh\big\}\times \big\{\cG_0=\ovl\cG\}$, which, by the discussion in Section~\ref{G0levels} and the choice~\equ{new domain} of the domain, are compact. Then, $\cM_{\hh, \ovl\cG,\varepsilon}$ remains compact for small $\varepsilon$, because collisions are excluded. Then, action--angles coordinates $\textrm{\sc a}=(\ovl\L, A, \ovl\l,\a)$ can be found in each connected component of $\ovl\cD$. In such coordinates, both $\hh_{\textrm{\sc k}}$ and $\cG_{\textrm{\sc k}}$  would depend on $(\ovl\L, A)$ only. Moreover, by its definition, $A$ is $\varepsilon$--close to
\beq{A0}A_0=\frac{1}{2\p}\int {\rm G}d{\rm g}\eeq
where ${\rm G}$ solves $\cG_0(\ovl\L,{\rm G},{\rm g})={\cal G}$. 
This function corresponds to $\ovl\L\sqrt{1-w^2}$, where 
has been computed in Section~\ref{G0levels}. From this expression, one sees that $\partial_{\ovl\cG}A_0\ne 0$ (being the integral of a  positive function , it is strictly increasing), therefore also $\partial_{\ovl\cG}A\ne 0$. By Implicit Function Theorem, one can invert $A$ as a function of $\ovl\L$, $\ovl\cG$. This allows to write $\hh_{\textrm{\sc a}}(\ovl\L,A(\ovl\L,\ovl\cG))=:\ovl\hh(\ovl\L,\ovl\cG)$, which corresponds to the sum of  series~\equ{eq: normalization}.

\subsection{$\ovl\l$--independence of $\ovl\cG$ and $(\ovl\UU$, $\cG_0)$ commutation}\label{lambda independence}
In this section we state an abstract result  that allows to prove (i) that $\ovl\UU$ commutes with $\cG_0$ and (ii) that $\ovl\cG$ is $\ovl\l$--independent.  Note that (i) easily implies~\equ{U}.

\nl
(i) and (ii) follow from the corresponding theses of lemma below, taking $\cH:=\ovl\hh$ in~\equ{eq: norm}, $\cJ:=\ovl\cG$ in~\equ{ovl G} and 
$(I, \f, p, q)=\ovl{\textrm{\sc k}}$.

\begin{lemma}\label{trivial lemma}
Let $(I, \f, p, q)$, with $(I,\f)$, $(p,q)$  pairwise conjugate, canonical coordinates on the phase space $\cP=V^1\times \torus^1\times \cB$, where $V^1\subset \real^1$, $\cB\subset \real^{2n}$ open and connected. Let $\cH:\ \cP\times (-\varepsilon_0,\varepsilon_0)\to \real$ a  $\f$--independent function of the form
 \beqa{close to be int}
\cH(I,p,q,\varepsilon)=\cH_0(I)+\varepsilon \cH_1(I,p,q)+\cdots
\eeqa
 analytic in $\varepsilon$. Then

\begin{itemize}
\item[{\rm (i)}]
 $\cJ_0$ Poisson--commutes with $\cH_1$;
 \item[{\rm (ii)}]  any  first integral  $\cJ(I,\f, p,q,\varepsilon)$ to $\cH(I,p,q,\varepsilon)$, analytic in $\varepsilon$,  with $\cJ_0(I,\f, p,q):=\cJ(I,\f, p,q,0)$ independent of $\f$, is $\f$--independent  for all $\varepsilon\in (-\varepsilon_0,\varepsilon_0)$.
 \end{itemize}
\end{lemma}
\proof
We start with (ii).
Let
\beqno\cJ(I,\f, p,q,\varepsilon)=\cJ_0(I, p,q)+\varepsilon \cJ_1(I,\f,p,q)+\cdots\ .\eeqno
We prove that  for all $j\ge 1$, $\cJ_j$ is $\f$--independent. By assumption, $\cI_0$ is $\f$--independent.  Assuming, inductively, that, for a given $k\ge 1$ and any $0\le j\le k-1$, $\cJ_j$ is $\f$--independent (so that the case $k=1$ corresponds with the assumption), we prove that $\cJ_k$ is so.
Writing the commutation relation of $\cH$ and $\cJ$, and picking the term proportional to $\varepsilon^k$ one as
$$\Big\{\cH_0,\cJ_k\Big\}+\sum_{j=1}^{k}\Big\{\cH_j,\ \cJ_{k-j}\Big\}=0$$
Since $\cH_0$ depends only on $I$, the first term in this identity has vanishing $\f$--average, while, by the inductive assumption, the second term is $\f$--independent. Then one has, identically,
$$\Big\{\cH_0,\cJ_k\Big\}=\sum_{j=1}^{k}\Big\{\cH_j,\ \cJ_{k-j}\Big\}=0$$
Therefore, $\cJ_k$ is $\f$--independent. The thesis (i)   follows from this identity, with 
$k=1$.

\subsection{The function $\widehat\UU$}\label{explicit}
In this appendix, we write relation~\equ{U} explicitly. This is not used in the paper, but, in view of equations~\equ{solutions} and Remark~\ref{rem}, may turn to be useful in applications.

\nl
As a matter of fact, there are infinite ways of representing $\ovl\UU$ as a function of $\rr'$, $\L$, $\Theta$, ${\cal G}_{0}$. Indeed, since $\ovl\UU(\rr',\L,\Theta,{\rm G},{\rm g})$ and $\cG_0(\rr',\L,\Theta,{\rm G},{\rm g})$ commute, if ${\rm G}(\cG_0, \rr',\L,\Theta,{\rm g})$ solves
$$\cG_0(\rr',\L,\Theta,{\rm G},{\rm g})=\cG_0$$
then the function $\ovl\UU\big(\rr',\L,\Theta,{\rm G}(\cG_0, \rr',\L,\Theta,{\rm g}),{\rm g}\big)$ is ${\rm g}$--independent, and one can take, for any fixed ${\rm g}_0$,
\beq{explicit form}\widehat{\rm U}\Big(\rr',\L,\Theta,{\cal G}_{0}\Big)=\ovl\UU\big(\rr',\L,\Theta,{\rm G}(\cG_0, \rr',\L,\Theta,{\rm g}_0),{\rm g}_0\big)\ .\eeq
A convenient choice is to take  ${\rm g}_0=\frac{\p}{2}$, in which case, as one sees from~\equ{G00}, ${\rm G}(\cG_0, \rr',\L,\Theta,{\rm g})\equiv\sqrt{\cG_0}$. Then~\equ{explicit form} becomes
$$\widehat{\rm U}\Big(\rr',\L,\Theta,{\cal G}_{0}\Big)=\ovl\UU\big(\rr',\L,\Theta,\sqrt{\cG_0},\frac{\p}{2}\big)\ .$$
More explicitly, using the expression of $\ovl\UU$ in terms of $\textrm{\sc k}$, which is
\beqano
\ovl\UU(\rr', \L, \Theta,{\rm G},{\rm g})&=&-\frac{1}{2\p}\int_{0}^{2\p}d\ell\Big[{{\rm r}'}^2\nonumber\\
&&+2 {\rr'} a\rr'\r \sqrt{1-\frac{\Theta^2}{{\rm G}^2}}\cos({\rm g}+\n)+a^2\rr^2\Big]^{-1/2}\eeqano
where $a$ is as in~\equ{a} and
\begin{itemize}
\item[--] $\ee=\sqrt{1-\frac{{\rm G}^2}{\L^2}}$ is the eccentricity;
\item[--] $\zeta$ is  the eccentric anomaly, solving Kepler's equation $\zeta-\ee\sin\zeta=\ell$;
\item[--] $\r:=|{\mathbf x}_{\textrm{\sc k}}|/a=1-\ee\cos\zeta$;
\item[--] $\n$ is the true anomaly, defined by $\n=\arg(\cos\zeta-\ee,\sqrt{1-\ee^2}\sin\zeta)$;
\end{itemize}
we find
    \beqa{ff}\widehat{\rm U}(\rr',\L,\Theta,{\cal G}_0)&:=&-\frac{1}{2\p}\int_{\torus}d\zeta(1-{\cal E}(\L,{\cal G}_0)\cos\zeta)\nonumber\\
    &&\Big[{\rr'}^2-2a{\rr'}{\cal I}(\L,\Theta,\cG_0)\sin\zeta\nonumber\\
    &&+a^2(1-{\cal E}(\L,{\cal G}_0)\cos\zeta)^2\Big]^{-1/2}\eeqa
where
\beqno {\cal E}(\L,{\cal G}_0):=\frac{\sqrt{\L^2-{\cal G}_0}}{\L}\qquad {\cal I}(\L,\Theta, {\cal G}_0):=\frac{\sqrt{{\cal G}_0-{\Theta^2}}}{\L}\ .\eeqno 
Note that negative values for the expressions under the square roots  are not  a problem, since $\widehat\UU$ is even in $\rr'$, $a$, $\cE$, $\cI$ separately.

\newpage\section{A Normal Form Theory}\label{normal form}
Let $\widehat\cD=\widehat\cB\times \torus^2$, with $\widehat\cB\subset \real^4$ compact and $\torus:=\real/(2\p\integer)$ be a real domain for the coordinates $\big((\widehat \rr',\widehat\L,\widehat A, \widehat\RR'),\widehat\ell,\widehat \a\big)$ where $\widehat\hh$ is regular, and let $\widehat\cD_\complex\supset\widehat\cD$ be a suitable complex, compact domain such that $\widehat\hh$  has an holomorphic extension on $\widehat\cD_\complex\supset\widehat\cD$.
Let us denote:
$$\o_{(\widehat\rr',\widehat\L,\widehat A)}:=\partial_{(\widehat\rr',\widehat\L,\widehat A)}\widehat\HH_0$$
$$ \langle \widehat f\rangle_{\widehat\ell,\widehat \a}:=\frac{1}{(2\p)^2}\int_{[0,2\p]^2}\widehat f d\widehat\ell d\widehat \a$$
$$\|g\|:=\sup_{\widehat\cD_\complex}|g|\ ,\ {\cal X}:=\sup_{\widehat\cD_\complex}|\widehat \RR'|\ .$$

\nl
The following result is known, even though not completely standard. The peculiarity of the Normal Form Lemma is that no small divisors condition is needed. This is possible because the coordinate $\widehat\RR'$ is not an angle. Roughly, this circumstance allows for a careful choice of the integration constant in the homological equation that allows to `de--singularize' the small denominators. Details will be published elsewhere.    In the text, it has been applied with $N\sim \varrho^{-1}$, as it is possible, since $\o_{(\widehat\L,\widehat A)}$ is of order $\varrho$, while $\o_{\widehat\rr'}$ is of order 1.

\nl
\begin{lemma}
There exists a constant $c$ such that, for any $N\in \natural$ such that the following inequalities are satisfied
\beqno c N{\cal X}\|\frac{\o_{\widehat\L}}{\o_{\widehat\rr'}}\|<
1
\ ,\
c N{\cal X}\|\frac{\o_{\widehat A}}{\o_{\widehat\rr'}}\|<
1\ ,\ c N{\cal X}   \|\frac{1}{\o_{\widehat\rr'}}\| \| \widehat f\|<1
 \eeqno
one can find a $\varrho$ --close to the identity canonical transformation
$$(\widetilde \rr',\widetilde\L,\widetilde A, \widetilde\RR',\widetilde\ell,\widetilde \a)\to (\widehat \rr',\widehat\L, \widehat A, \widehat\RR',\widehat\ell,\widehat \a)$$
 that carries $\widehat\HH$ to
$$\widetilde\HH=\widehat\HH_0+\widetilde\HH_1+\widetilde f$$
where $\widetilde\HH_1$, $\widetilde f$ satisfy
\beqa{thesis}
\|\widetilde\HH_1-\langle \widehat f\rangle\|\le c {\cal X}\|\frac{1}{\o_{\rr'}}\| \| \widehat f\|^2\ ,\     \|\widetilde f\|\le \frac{\varrho^2}{2^{N+1}} \|\widehat f\|\ .\eeqa
\end{lemma}

\nl The formula in~\equ{2body} follows from the former inequality in~\equ{thesis} and $\langle{\mathbf y}'_{\textrm{\sc a}}\cdot{\mathbf y}_{\textrm{\sc a}}\rangle_{\widehat\ell,\widehat\a}=0$, so that
\beqa{pert}\langle \widehat f\rangle_{\widehat\ell,\widehat\a}=\frac{\langle|{\mathbf y}'_{\textrm{\sc a}}|^2\rangle_{\widehat\ell}}{2\mm'{\rr'}^2}=\frac{\widehat{\RR'}^2}{2\mm'{\rr'}^2}+\frac{\langle(\GG-{\rm G})^2\circ\textrm{\sc a}\rangle_{\widehat\a}}{2\mm'{\rr'}^2}\eeqa

\newpage
\addcontentsline{toc}{section}{References}
\def\cprime{$'$} \def\cprime{$'$} \def\cprime{$'$} \def\cprime{$'$}

\end{document}